\documentclass[runningheads]{llncs}

\usepackage{graphicx}        % standard LaTeX graphics tool
\usepackage{algorithm2e}
\usepackage{colortbl}
\usepackage{color}
\definecolor{Red2}{rgb}{1,0.7,0.7}
\definecolor{Green2}{rgb}{0.7,1,0.7}
\definecolor{Blue2}{rgb}{0.85,0.85,1}
\definecolor{Yellow2}{rgb}{1,1,0.4}
\newcolumntype{?}{!{\vrule width 0.8pt}}
\makeatletter
\def\hlinewd#1{%
\noalign{\ifnum0=`}\fi\hrule \@height #1 %
\futurelet\reserved@a\@xhline}
\makeatother

\begin{document}
\title{Scattered data approximation by LR B-spline surfaces. A study on refinement strategies for
efficient approximation}
\titlerunning{Refinement strategies for LR B-spline approximation}
\author{Vibeke Skytt \and Tor Dokken}
\institute{SINTEF, Forskningsveien 1, P.O. Box 124 Blindern, 0314 Oslo, Norway}
\maketitle
%\address{SINTEF, Forskningsveien 1, 0314 Oslo, Norway}

\begin{abstract}
  Locally refined spline surfaces (LRB) is a representation  well suited for  scattered data approximation. When a data set has local details in some areas and is largely
  smooth in other, LR B-splines allow the spatial distribution of degrees of freedom to follow the variations of the
  data set. An LRB surface approximating a data set is refined in areas where the accuracy does not meet
  a required tolerance. In this paper we address, in a systematic study, different LRB refinement strategies and polynomial degrees for surface approximation. We study their influence on data volume and accuracy when applied to geospatial data sets with different structural behaviour. The performance of the refinement strategies is to some degree coherent and the article concludes with some recommendations.
  An overall evaluation indicates that bi-quadratic LRB are preferable for the uses cases tested,  and that the strategies we denote as "full span" have the overall best performance.
\end{abstract}

\section{Introduction}
Tensor-product B-spline surfaces is a mature and standardized geometry representation that has been known at least since the 1970s. The first uses of Tensor-product B-spline  were in Computer Aided Design (CAD). In Isogeometric 
Analysis (IgA)~\cite{IGA} B-splines replace the traditional shape functions used in Finite Element
Analysis (FEA).
A guide to  splines can be
found in~\cite{splines}. Tensor-product spline surfaces have very good
numerical properties, but lack  local refinement of the spline space. In recent years several approaches have been proposed for local refinement of spline surfaces including T-splines~\cite{Tsplines2,Tsplines1},
Truncated Hierarchical B-splines (THB)~\cite{hierarchical1,hierarchical2}, and LR B-splines~\cite{LRbasic}.
We will in this paper focus on the use of LR B-spline surfaces for approximation of scattered data.

We perform a systematic study on the effects of different strategies for local refinement of LRB surfaces for the approximation of geospatial point clouds. We focus on LRB as LRB allows more flexibility of refinement strategies than T-splines and THB.
Section~\ref{sec:local} presents locally refined splines in general and LRB in particular.
We then turn to scattered data approximation in Section~\ref{sec:approx}. Section~\ref{sec:refstrat} gives a brief
overview over previously published local refinement strategies for bi-variate splines and
discuss the
concept of a good refinement strategy for approximation of large point clouds. A set of candidate
strategies are defined in~\ref{sec:candstrat}. Section~\ref{sec:study} presents five
data sets and the corresponding approximation results for the selected refinement strategies along with
analysis of the result related to each data set. 
In Section~\ref{sec:summary} threads from the case
specific analysis are drawn to get a unified understanding and finally, in Section~\ref{sec:conclude} some
conclusions are drawn.

\subsection{Locally refined splines} \label{sec:local}
The lack of local refinement for tensor-product splines provides  severe restrictions for IgA as well as for scattered data approximation. Using tensor product B-splines will in most cases introduce significantly more  degrees of freedom than actually needed. This makes the data volume grow and restricts the size of problems to be addressed. Two basic approaches exist  local spline:
\begin{itemize}
    \item Refinement in the mesh of vertices/control points, the approach for T-splines.
    \item Refinement in the parameter domain, the approach for THB
    and LRB.
\end{itemize}
It is attractive to span the locally refined spline space by a basis,  a set of functions that are linearly independent. 
T-splines, LRB and THB are all spanned by functions that are composed from tensor products of univariate B-splines.

\subsubsection{T-splines.}\label{sec:T-splines}
T-splines were introduced by Sederberg et. al.~\cite{Tsplines1} to enhance the modelling flexibility in CAD design. New control points are added to the what is denoted a T-mesh. The starting point is a tensor product B-splines  surface and its regular grid of control points. Then additional control points are added according to a set of T-spline rules. The knot vectors of the tensor product B-spline corresponding to a control point are defined by traversing the T-mesh  starting from the control point and going outwards in all four axes parallel directions until a degree dependent number of meshlines are intersected. In the bi-cubic case this traversing stops after two lines in the T-mesh are intersected in each of the four directions.    The new control points are used to model local details in a preexisting surface. Frequently the B-splines spanning the T-spline space have to be scaled to form a partition of unity.

The most general version of T-splines possesses neither nested spline spaces nor a guarantee for linear
independence~\cite{T_lindep} of B-splines. In IGA, linear independence is important and
additional rules were added to the T-spline creation algorithm to ensure
Analysis Suitable T-splines~\cite{Tsplines2}.

\subsubsection{Truncated Hierarchical B-splines.}\label{sec:THB}
Hierarchical B-splines (HB), introduced by Forsey and Bartels \cite{hierarchical1}, are  based on a dyadic sequence of grids determined by scaled lattices over which uniform spline spaces are defined. HB provides nested spline spaces,  is spanned by tensor product B-splines but do not provide spanning functions that are a partition of unity. How to select B-splines that gives HB linear independence was solved in \cite{HBindependence}.  To provide partition of unity of HB the Truncated Hierarchical B-splines (THB)~\cite{hierarchical2} were introduced. THB are linearly independent and reduce the footprint of the basis functions. The basis functions of THB are made by truncating tensor product B-splines with tensor product B-splines from a finer refinement levels. Some B-splines are not truncated and remain tensor product B-splines. The truncated B-splines can be described as a sum of scaled tensor product B-splines from finer levels (or alternatively as a difference between the tensor product B-spline truncated and scaled tensor product B-splines from finer levels).

\subsubsection{Locally Refined B-splines.}\label{sec:LRB}

An LR B-spline surface is a piecewise polynomial or piecewise rational polynomial surface defined on an LR-mesh. The
LR-mesh is limited by the configurations achievable by applying a sequence of
refinements starting from a tensor-product mesh. The LR B-spline representation is, consequently, algorithmically defined.
An LR B-spline surface is defined as
$$
F(u,v) = \sum_{i=1}^L P_i s_i R_{i,p_1,p_2}(u,v) \, ,
$$
where $P_i$, $i=1,\ldots ,L$ are the surface coefficients, and $s_i$ are scaling factors introduced
to ensure partition of unity of the resulting collection of tensor product B-splines.
The tensor product B-splines $R_i$ are  of bi-degree $(p_1,p_2)$ defined on
knots extracted from the knot vectors in the $u$ and $v$ direction of the parameter domain.

\begin{figure}
\begin{center}
\includegraphics[width=0.5\textwidth]{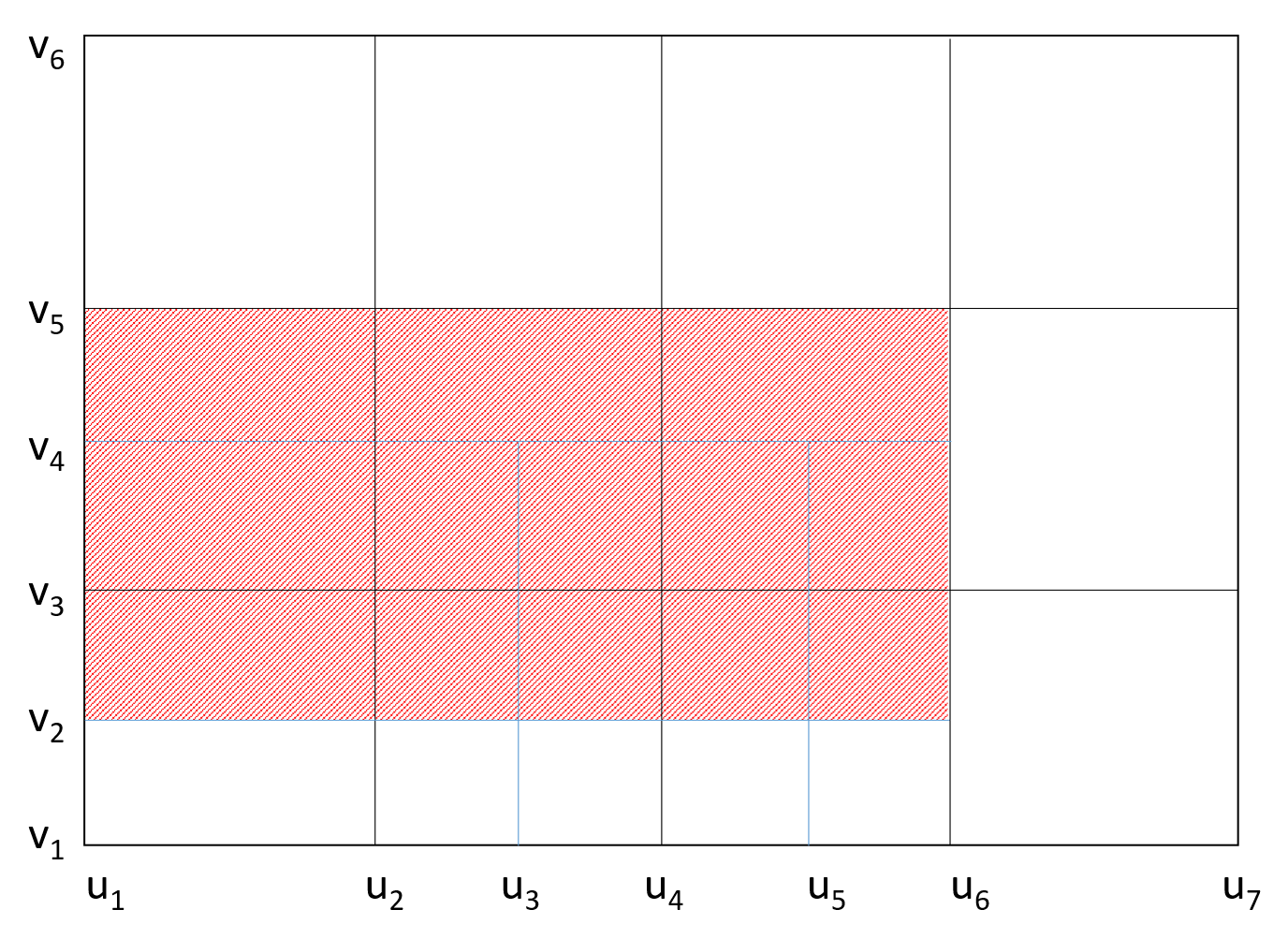}
\end{center}
\caption{\label{fig:lrmesh} LR-mesh. The support of one tensor product B-spline visualized as a red pattern.
  Initial knotlines are shown as black lines, the inserted knotline segments are blue.}
\end{figure}
The initial tensor-product mesh corresponding to the LR-mesh showed in Figure~\ref{fig:lrmesh} is
given by the knot vectors 
$[u_1, u_1, u_1, u_2, u_4, u_6, u_7, u_7, u_7]$ and \newline
$[v_1, v_1, v_1, v_3, v_5, v_6, v_6, v_6]$. 
It corresponds to a polynomial spline surface of degree 2 in both
parameter directions with multiple knots in the end parameters. The LR-mesh is constructed by first
inserting knots at $v_2$ and $v_4$ covering a part of the surface domain, and next inserting
knots at $u_3$ and $u_5$. The B-spline with support shown in red has the knot mesh
$[u_1,u_2,u_4,u_6] \times [v_2, v_3, v_4, v_5]$. Some knotlines in the support of the B-spline are not a part of the B-spline
definition as they do not traverse the complete support of the tensor product B-spline.

The procedure for refining an LR B-spline surface is as follows:
\begin{enumerate}
\item Select a new knotline segment traversing the support of at least one tensor product B-spline. The segment can be an
  extension to an existing knotline segment. Thus, a knotline going from $(u_5,v_4)$ to $(u_5,v_4)$ is a legal
  choice.
\item Subdivide all tensor product B-splines where the new knotline/extended knotline constitutes a legal inner knotline segment.
\item The tensor product B-splines are required to have minimal support. After
  inserting a new knotline segment and performing the subdivision in 2., there might still be tensor product B-splines that do not have minimal support. Consequently all such B-splines  are split accordingly.
  This process is continued until all tensor product B-splines have minimal support.
\end{enumerate}
LR B-spline surface refinement is described in detail in~\cite{IGAref1}.
If more than one new knotline segment is defined simultaneously, the refinement process is applied to
one segment at the time.

Nested spline spaces are ensured by the LR B-spline construction. The LR B-splines are non-negative and have
compact support. The scaling factors $s_i$ are computed during the refinement process to ensure partition of unity. 
The LR B-splines are
not guaranteed to be linearly independent, but a dependency relation can be detected and resolved by
dedicated knot insertions. Linear dependency can only occur in a situation with overloaded LR B-splines.
An element is overloaded if it belongs to the support of more LR B-splines than necessary to span the
spline space over the element. An LR B-spline is overloaded if all elements in its support are overloaded.
Overloading can be detected by the peeling algorithm~\cite{LRbasic}.
Patrizi and Dokken address configurations that can
lead to a linear dependency relation in~\cite{LinDep}. 
In a situation with local linear independence over all polynomial elements all weights, $s_i$ $i=1, \ldots L$, will be equal to one.

\subsection{Scattered data approximation} \label{sec:approx}
The aim is to approximate a scattered data point cloud by an LR B-spline surface. We focus on
data sets with projectable points that can be parameterized by their $xy$-coordinates leaving
the $z$-coordinate to be approximated by a height function. Also non-projectable points being
parameterized by some appropriate method can be handled by the general algorithm used in this paper.

\begin{algorithm}
\KwData {Point cloud, maximum number of iterations, tolerance}
\KwResult{Approximating surface, information on approximation accuracy}
Generate initial surface\;
Compute accuracy\;
\While{there exists points with larger distance than the given tolerance,
the maximum number of iterations is not reached and new knotlines can be inserted}
{
  Refine the surface where needed\;
  Perform approximation in the current spline space\;
  Compute accuracy\;
}
\caption{Iterative algorithm for LR B-spline surface generation \label{alg:approx}}
\end{algorithm}
Algorithm~\ref{alg:approx} gives an overview of the iterative approximation algorithm. The starting
point is a tensor-product B-spline surface defining an initial spline space. It is beneficial if the
corresponding coefficients give a rough indication on the behaviour of the shape represented by
the point cloud, but is is not a requirement for the algorithm.

The focus of
this article is to investigate the term ``Refine the surface where needed''. The polynomial segments of the surface will, during
the computation, keep track of the data points situated within their domain as well as  some accuracy
statistics. This includes maximum and average distances between the surface and the points in  the segment,
and the number of points with a distance to the surface that is larger than the specified tolerance. This information
provides a basis for selecting where and how new degrees of freedom shall be added. In the following  a polynomial segment will be denoted an element. 

Two surface approximation methods are applied during the iterative algorithm: Least
squared approximation with a smoothing term and multi-level B-spline approximation (MBA).
Let ${\mathbf x}=(x_k,y_k,z_k), k=1,\ldots,K$ be the projectable point cloud
we want to approximate. Least squares approximation is a global method where a minimization
functional
$$
\min_{\mathbf c} [\alpha _1 J(F(x,y) + \alpha _2 \sum _{k=1}^K (F(x_k,y_k) - z_k)^2]
$$
is differentiated and lead to a sparse, linear equation system. Linear independence of the B-splines
is a prerequisite for a non-singular equation system. $F(x,y)$ denotes the
LR B-spline height function we want to obtain. The solution to this functional
results in the best possible
approximation with the given degrees of freedom in a least squares sense. The actual
smoothing term $J(F(x,y)$ is described in~\cite{lralg1} and other possible smoothing terms can be found
in~\cite{smoothterm}. Our focus is on the approximation. As the input points may or may not
represent a smooth surface, the weight $\alpha _1$ on smoothing must be kept low. Still the term is
important to handle areas in the surface domain without points.

The multi-level B-spline approximation algorithm were introduced by Lee et. al. in~\cite{mba} 
for scattered data interpolation. It is an iterative, local approximation method that results in a
hierarchical structure of tensor-product B-spline surfaces. At each iteration level the residuals between
the point cloud and the previous surface is computed using a pseudo inverse approach.  If 
${\mathbf x_c}=(x_c,y_c,z_c), c=1,\ldots,C$ are the residuals situated
in the domain of one B-spline, the coefficient of the residual surface $P_i$ is determined as
\begin{equation}
P_i = \frac{\sum_{c=1}^C N_c(x_c,y_c)^2 \phi_c}{\sum _{c=1}^C N_i(x_c,y_c)^2}, 
\label{eq:MBA}
\end{equation}
where
$$
\phi_c = \frac{N_i(x_c,y_c)z_c}{\sum_{j=1}^J N_j(x_c,y_c)^2} \, .
$$
$N_j, j=1,\ldots,J$ are the B-splines overlapping the residual $(x_c,y_c,z_c)$.
$\phi_c$ is obtained by a pseudo inverse approach
to solve the under-determined equation system 
$z_c=\sum_{j=1}^J \phi_j N_j(x_c,y_c)$. Every residual lead to different
values of $\phi$ and the expression (\ref{eq:MBA}) is obtained by minimizing
the error $e(P_i) = \sum_{c=1}^C |P_i N_i (x_c,y_c)-\phi_i N_i (x_c,y_c)|$
with respect to $P_i$.
The process is explained in more detail in~\cite{mba2}.

In the LR B-spline setting, the residual surface is incorporated in the expression for the approximating
surface avoiding a hierarchical representation.
The MBA method does not require linearly independent basis function to converge to a solution.

Both approximation methods are in some sense minimizing the average distance between the point cloud
and the surface. Thus, neither the maximum distance nor the number of points outside the tolerance
can be expected to decrease constantly. The average distance will in general be steadily reduced, but
temporary stagnation may occur in particular in the context of outliers or if the elements have been
refined to such a degree that we model noise.
In general we expect the input point clouds to have high data sizes and a varying degree of smoothness
over the domain of the point cloud, which makes the property of local refinement essential.

\section{Refinement strategies and success criteria} \label{sec:refstrat}
Locally refined spline surfaces are appropriate for approximation of large scattered data sets
due to the ability of increasing the data size in the areas where more degrees of freedom is required while
keeping the data size low in other areas. However, it is not obvious how to decide where 
new knotline segments should be inserted. There is a large degree of freedom in selecting
these knotlines, and we will investigate the effect of different choices.

Previous studies on refinement strategies for LR B-splines have focused on the use of
LR B-splines in isogeometric analysis or refinement strategies that ensure local or global linear
independence. Local linear independence ensures that a minimum amount of B-splines overlap an element
while
global linear independence implies that a linear equation system originating from a least
squares approximation or finite element analysis is non-singular.

Johannessen et. al. give an introduction to LR B-splines in the context of isogeometric analysis in~\cite{IGAref1}
where a detailed description of the refinement basics is provided and the effect on a number of test cases is investigated for the
three  following refinement strategies:
\begin{description}
  \item [Full span:] Given an element selected for refinement, all B-splines overlapping this element
    are split in a parameter contained in the element.
  \item [Minimum span:] The shortest possible knotline overlapping a chosen element that splits
    at least one B-spline, is selected. Several candidate B-splines can exist and even if additional
    selection criteria are added, the candidate may not be unique. The selected knotline is not necessarily
    symmetric with respect to the element.
  \item [Structured mesh:] Choose a B-spline and refine all knot intervals in this B-spline.
    No elements with an unbalanced aspect ratio will be created if the split is performed in the
    middle of the knot intervals.
\end{description}
The sensitivity towards the different choices of refinement strategies were in \cite{IGAref1} found to be
moderate, but high regularity meshes tended to give less error in the computation compared to lower
regularity meshes for the same number of degrees of freedom.
In~\cite{IGAref3} the structured mesh approach is analyzed theoretically and
numerically in a set of test cases.

Bressan and J\"uttler~\cite{Linref1} look at refinement of LR B-spline surfaces from the perspective of
local linear independence and present a mesh construction where this property is proved. Patrizi et.
al.~\cite{Linref2} propose a practical refinement strategy where local linear independence is
ensured. A modified structured mesh refinement where some of the splits are prolonged so the refined
mesh satisfies the so called non-nested support property is proposed, see~\cite{Linref2} for a definition.

Other authors have focused on refinement of hierachical B-spline surfaces or T-spline surfaces.
Bracco et. el. define two classes of admissible meshes for hierarchical B-splines and 
compare them for use in isogeometric analysis in~\cite{IGAref2}. A comparison between two refinement
strategies applied to hierarchical B-splines and T-splines in the
context of IGA is investigated in~\cite{IGAref4}. 

The applications of scattered data approximation and isogeometric analysis have some fundamental
differences. In the context of IGA, the refined meshes normally belong to an intermediate stage in
the computations. Thus, the possibly large locally refined spline surfaces are not kept. Moreover, the
need for refinement in an adaptive isogeometric analysis computation is typically concentrated in
localized areas. One of the main motivations for approximating scattered data with a locally refined
surface is the need for a compact representation of data with a locally non-smooth behaviour. Areas where
extra degrees of freedom are required for an accurate approximation, may be scattered around in the
entire domain of the data set.

We restrict ourselves to polynomial LR B-spline surfaces of bi-degree one, two and three and place all
new knotline
segments in the middle of an existing knot intervals. With this configuration, no linear dependence
relation has been encountered.

If there exist data points within an element with a distance to the surface larger than the
tolerance, there is either a need for more degrees of freedom corresponding to the element,  or
the required accuracy cannot be met by a smooth surface due to outliers or a lack of smoothness
in the data points.
Any B-spline overlapping this element will give new degrees of freedom to the
element, if split. This implies that there is freedom in how we choose to increase the degrees of
freedom. It is also clear that all choices are not equally good.

We define some criteria for a good refinement strategy:
\begin{enumerate}
\item Best possible accuracy with the minimum degrees of freedom
\item If supported by the data set, it should be possible to adapt the surface to
  the point cloud within the prescribed tolerance
\item Avoid a premature stop
\item A steady improvement in accuracy with new degrees of freedom
\item Keep the execution time low \label{time}
\item Keep the memory consumption low \label{memory}
\item A new knotline in the middle of an element, should not to a large extent lead to adjacent elements being split at other positions than in the middle 
\item The resulting surface has a distribution of knotlines reflecting the behaviour of the data set
  and the difference in B-spline size is minimized. B-splines do not contain
  many unused/not needed knotlines. \label{Bspline_prop}
\item Linear independence and local linear independence of B-splines
\end{enumerate}
We will mainly focus on the first five criteria. The Critera~\ref{time} and~\ref{memory}
are to some degree linked and also dependent on the
previous ones. A lean surface will lead to lower memory consumption and the part of the
execution time spent in surface refinement is linked to surface size. The execution time
also depends on the number of steps applied in the iterative algorithm and whether or not
the knot insertion at late iteration steps is focused in a few areas or spread out in
the entire surface domain.

The importance of linear independence depends on the use of the
resulting surface. The approximation algorithm outlined in Algorithm~\ref{alg:approx} is not dependent
on linear independence. We will, thus, not focus on the last criterion.

We can now formulate two hypothesises on properties of
a good refinement strategy. They are taken into account
when the strategies to be tested are defined in the
next section:
\begin{description}
\item [Hypothesis 1:] A gradual introduction of new degrees of freedom gives better
approximation efficiency.
\item [Hypothesis 2:] An improved accuracy can be blocked by failing to split one or 
more important B-splines.
\end{description}
An important term in this discussion is {\bf approximation efficiency}. It is defined as the number
of resolved points divided by the number of surface coefficients for a particular refinement
strategy and a particular iteration level. This figure will, together with other criteria, be
used to evaluate the success of one refinement strategy compared to others. 

\section{Applied refinement strategies}\label{sec:candstrat}
\subsection{Main categories}
The selected data sets are approximated using a variety of refinement strategies. Refinement can
be triggered either from an {\bf element} or a {\bf B-spline} containing unresolved data points
in its domain. A new knotline segment must split at least one tensor product B-spline (from now on denoted a B-spline).

For strategies triggered by an element,
one or more B-splines overlapping this element must be selected. If all
B-splines overlapping the element are selected, the corresponding strategy is named a {\bf Full
  span} strategy.

A {\bf Minimum span} strategy selects one overlapping B-spline based on some criteria.
The new knotline splits the element in half in the direction of refinement  and it is
dimensioned to split all selected B-splines.
We will test three different minimum span strategies. They differ in how the  B-spline
to split is selected. The selection criteria are: the {\bf largest} overlapping
B-spline, the overlapping B-spline with the highest number of {\bf unresolved} points in the
B-spline support compared to the total number of points in the support, and a {\bf combination}
of the two where the criteria are given equal weights. Several B-splines may fit the
selected criterion equally well. In case of doubt, the most centered B-spline with respect to the
initial element is selected. If there is still no unique B-spline satisfying the criterion, one
candidate is randomly chosen.

A refinement strategy triggered by a B-spline can impose refinement in all knot intervals in the
B-spline. This is denoted a {\bf Structured mesh} strategy, and is similar to the refinement strategy for hierarchical B-splines. Alternatively, the refinement can be
performed only in knot intervals in the support of the B-spline where the corresponding elements have a significant number of points
outside the tolerance  or where the distances between the point cloud and the surface is
large. This is denoted a {\bf Restricted mesh} strategy. 

\subsection{Restrictions to the introduction of new knots}
We want to test the hypothesis that a gradual introduction of new degrees of freedom will lead to
a lean final surface with a good accuracy.
The possible number of new knots introduced at each iteration level can be reduced by applying
refinement in {\bf alternating parameter directions}. The surface is refined in the first parameter
direction at odd levels and in the second parameter direction at even iteration levels. Refinement
can also be performed at {\bf both parameter directions} at each level.

{\bf A threshold} can also be applied to restrict the number of new knotlines at each step,
possibly in addition to alternating parameter directions and/or restrictions inherent in the
strategy itself. Thresholds can imply that only the elements or B-splines with the most
significant approximation errors, either with respect to distance or number of unresolved points, is
selected for refinement. The threshold factor is set globally for each iteration level.
In B-spline based strategies can thresholds be used to reduce the number of knot intervals that
are refined in
the selected B-spline. A refinement strategy can be combined with zero, one or two different
types of thresholds.

\subsection{Modifications to particularly restricted strategies}
The restricted mesh strategy may fail to capture elements with significant
points outside the tolerance belt where this element is atypical in the overlapping B-splines selected for splitting.
This situation can imply that
refinement in the associated knot interval is not performed and occurs typically at the
border of the point set. To ensure refinement in such elements, an element extension to
the strategy is applied. Elements with a significant approximation error that are not already split, are identified at each iteration level and trigger an additional 
element based refinement step. This approach is related to the full span strategy,
but the significance
of each parameter direction is evaluated individually. Thus, a strategy to refine in
both parameter directions will not necessarily imply that exactly these elements trigger refinement
in both directions. 

In addition to the restricted mesh strategy  the various minimum span strategies can lead to few new
knotlines being inserted at each iteration level, in particular if the strategy is combined with a
threshold.
We investigate the effect of combining these strategies with a full span strategy as the full span
strategies are unambiguous and lead to longer new knotline segments than the minimum span and the
restricted mesh strategies. The two combined strategies are applied with refinement in the same
number of parameter directions. A threshold with respect to the number of points outside the
tolerance in an element or B-spline domain is applied.  The actual type depends on the current strategy. When the convergence of the
initially chosen strategy slows down  the algorithm will switch to the corresponding full span
strategy.
The switch is performed when the fraction between newly resolved points in the last iteration and
the number of new coefficients drops below 0.1. 

\subsection{The components for refinement strategies}
Each refinement strategy and each combination of strategy and threshold is given a unique
label. The label is composed by one letter for each component in the strategy including threshold.
The components are:
\begin{description}
\item [Trigger entity] Element (e) or B-spline (b)
\item [Element based strategy] Full span (F) or Minimum span (M)
\item [B-spline based strategy] Structured mesh (S) or Restricted mesh (R)
\item [Element extension to B-spline based strategy] Select significant elements left out in the B-spline
  refinement (L)
\item [Minimum span selection criterion] Largest overlapping B-spline (l), overlapping B-spline with
  most unresolved points (u) or a combination (c)
\item [Parameter direction] Refine in both (B) or alternating (A) parameter directions at each iteration level
\item [Threshold] With respect to distance (td), the number of unresolved points in an element (tn)
  or the number of unresolved points in all elements corresponding to a B-spline knot interval (tk)
\end{description}

The size of the threshold type td 
is defined as a weighted sum between the maximum distance between the complete point cloud
and the surface, the average distance in unresolved points and the tolerance. An element or B-spline
is selected for refinement if the maximum distance between the points in the support of this
entity and the surface exceeds this threshold.

The type tn threshold is defined from the number of
unresolved points in an element combined with the distance between the surface and these points.
The minimum and maximum value over all elements are computed and the threshold is a weighted sum of
these two numbers. Elements with a higher score than the threshold trigger refinement.

The type tk
threshold is applied to the restricted B-spline based strategies only.
The threshold factor is set to the average
number of out-of-tolerance points in all B-spline supports
unless the average number of points in these supports is high.
Then the threshold factor equals 1\% of the average number of points in the supports.
Knot intervals in a selected B-spline is refined only if the number of out-of-tolerance points
in the elements belonging to that interval exceeds this threshold factor. This is a very strong threshold
as the number of points outside the tolerance
in the elements belonging to one knot interval in one of the parameter directions
is compared to the number of points in the entire B-spline support. In general are the threshold types and associated
factors set experimentally and the degree of tuning is low.
The aim is to get an impression on how different approaches influence the approximation process. 
All threshold factors are reduced with successive iterations.

\begin{table}
%\setlength{\tabcolsep}{4pt}
%%\centering
%  \footnotesize
\caption{Composition of strategy labels.}
\label{tab:strategies}
\begin{tabular} {|l|l|l|l|l|l|p{2cm}|l|p{2.7cm}|}
  \hline
  base & type & sub type & b2 & t2 & dir & threshold & e.g. & comment \\
  \hline
  e & F & & & & A or B & none, td or tn & eFA tn & Full span \\
  e & M & l or u or c & & & A or B & none, td or tn & eMlB & Minimum span \\
  e & M & l or c & (e) & F & A or B & tn & eMc/FA tn & Combine minimum and full span \\
  \hline
  b & S & & & & A or B & none or td & bSB td & Structured mesh \\
  b & R & & & & A or B & none, td, tk or td and tk & bRA td+k & Reduced mesh \\
  b & R & & e & F & A or B & tk and tn & bR/eFB tk/n & Combine reduced mesh and full span \\
  b & R & & e & L & A or B & none, td, tk or td and tk & bR+eLA tk & Element extended reduced mesh \\
  \hline
\end{tabular}
\end{table}
\subsection{Composition of refinement strategies for the tests}
The refinement strategies are applied to a number of test cases. Each strategy is tested with and
without threshold and with different choices in the
number of parameter directions to refine at each 
iteration step. Table~\ref{tab:strategies} shows how the strategy labels are composed
from its basic components. Here the base refers to whether the strategy is element or B-spline based and
the type further specifies the strategy to be full span or minimum span in the element case, or
structured or reduced in the B-spline case. The sub type applies to the minimum span strategy only.
If several strategies are applied during the iterations, a second combination of base and type
(b2 and t2) is specified.
The two strategies can be applied simultaneously as is the case for
the element extended restricted mesh strategy, or successive. If two element based strategies are
applied, the second e is omitted in the label. Two strategies or thresholds that are applied to
the same iteration level is combined with + in the label, and with / if they are used at different
levels. An extensive set of combinations of strategies and threshold types
are tested, but there are still possible combinations that are not applied. We believe, however,
that the results give a broad basis for drawing conclusions on the suitability of the various
strategies. All strategies are tested with alternating strategies and with refinement in both
parameter directions at every iteration level.

\section{Refinement study} \label{sec:study}
Five geospatial data sets are  used in the study. All data sets have points that are projectable on the $xy$-plane,
but have otherwise different sizes and properties.
The number of points vary from 71\,888 to
131\,160\,220. The regularity distribution varies from completely regular to scanlines were the
distances between scans are very large compared to the distances between points in the scans.
The smoothness also varies from data set to data set. Two of the data sets contain known outlier points.
The point clouds are approximated using various combinations
of refinement strategies and thresholds. All other factors in the computations are kept
constant. All data set, but one, is tested with a tolerance of 0.5 meters. The number of iterations
varies according to properties of the data set, and to some extent with regard to whether
alternating parameter directions are applied. Varying the tolerance and the number of iterations for
each data set could bring additional information to the range of refinement strategies, but fall
outside the scope of this article. The study is run using bi-degrees (1,1), (2,2) and (3,3).

The information related to each data set includes measurements of execution time.
The recorded time includes all aspects
of the approximation, but excludes reading and writing to file.
The computations are performed on a stationary desktop with 64 GB of DDR4-2666 RAM. It has a i9-9900K CPU with 8 cores and 16 threads, but a single core implementation is used in the 
experiments.
Approximation results are collected at the final stage of the
computations and at an intermediate stage. This stage will be defined according
to properties of each data set. 

Tabulated results on the accuracy of all refinement strategies,
combination with threshold types
and polynomial bi-degrees one, two and three will be made available in a separate report along with more
analysis of the results. 

\subsection{Banc du Four}
\subsubsection{Data set}
\begin{figure}
  \centering
\begin{tabular}{ccc}
\includegraphics[width=0.32\textwidth]{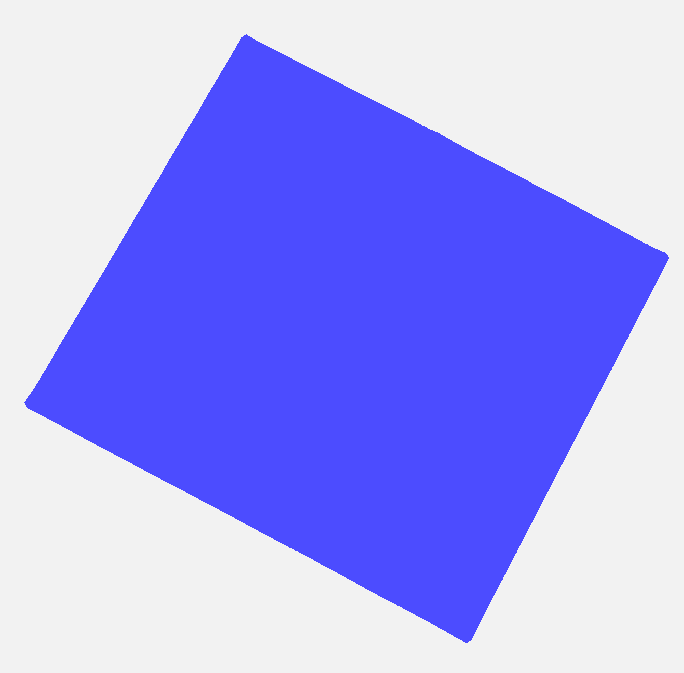}
&\includegraphics[width=0.3\textwidth]{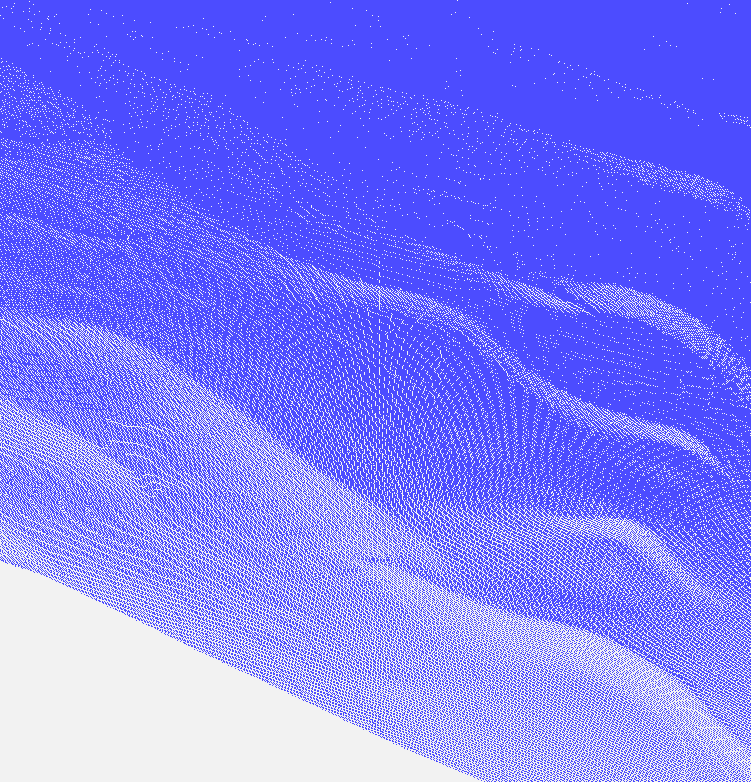}
&\includegraphics[width=0.32\textwidth]{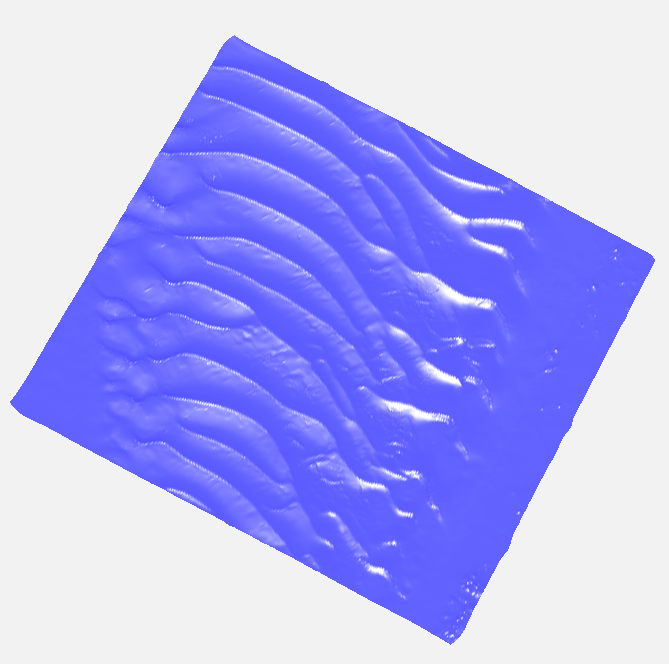} \\
(a) & (b) & (c)
\end{tabular}
\caption{(a) Sand dunes, France. (b) A detail. (c) Corresponding LR B-spline
  surface, the surface is trimmed to adapt to the point cloud domain. \label{fig:Fourpts}}
\end{figure}
Figure~\ref{fig:Fourpts} (a) shows a sub sea data set representing sand dunes, (b) zooms in on a detail.  The data
is obtained from the Banc de Four outside the coast of Brittany. The data set consists of 5\,054\,827
well distributed points. The elevation range is $[-92.76, -54.84]$ and the standard
deviation is 8.08. Sand dunes has a relatively smooth shape and the point cloud does
not contain any outliers. The approximating surface, Figure~\ref{fig:Fourpts} (c),
has an accuracy of 0.2 meters. 

\subsubsection{Result overview}
We use a tolerance of 0.2 meters for this data set.
The approximation algorithm is allowed to run for 40 iterations. It stops either because
the requested accuracy is reached, the maximum number of iterations is met or no suitable
new knotlines are identified. 

\begin{table}
  \caption{Computation time and number of coefficients. Performance
  ranges and associated strategies for polynomial degrees one to three. }
  \label{tab:F_range}
  \begin{tabular}{|l||l|l||l|l||l|l||l|l|}
    \hline
Deg & Min t. & Strategy & Max t. & Strategy & Min cf & Strategy & Max cf & Strategy \\
\hline
\hline
1 & 0m21s & bSB & 1m32s & eMuA td & 42\,503 & bR+eLA tk & 91\,050 & bSB \\
\hline
2 & 0m33s & eFB, bSB & 2m58s & eMuA td & 39\,366 & bR+eLA tk & 111\,992 & bSB \\
\hline
3 & 1m6s & eFB & 5m46s & eMuB td & 45\,826 & bRB tk & 162\,972 & bSB \\
\hline
\end{tabular}
\end{table}

The data set is well adapted for approximation with LR B-spline surfaces and most combination
of refinement strategies and threshold types
converge within the prescribed number of iterations. The execution time and
final number of coefficients, however, differ. Table~\ref{tab:F_range}
presents the minimum and maximum figures for execution time and 
the final number of coefficients along with the associated refinement strategies. The structured mesh strategy appears both as the 
strategy with the lowest execution time and the highest number of
coefficients. The minimum number of coefficients are found for
restricted mesh strategies while a minimum span strategy with
threshold based on distance has the highest computational time.

Several varieties of the minimum span refinement
strategies struggle  with convergence in the bi-linear case. Full span
refinement with a threshold also fails to reach the requested accuracy for some combinations. Most B-spline based strategies converge, the
only exceptions are restricted mesh with alternating parameter directions
and threshold involving the number of outside points (bRA td and bRA td+k). 

Almost all refinement strategies converge in the quadratic case. The exception is minimum span strategies where the overlapping B-splines are selected from
the number of out-of-tolerance points in the element (eMu) and the same restricted B-spline based
strategies as in the bi-linear case addressed above.

Bi-cubic polynomials result in more coefficients and higher execution times than lower
order polynomials, but most refinement strategies converge. The exceptions are also in this
case minimum span strategies with subtype u and restricted B-spline based strategies with threshold.

\subsubsection{Best performances}
\begin{table}
%\setlength{\tabcolsep}{4pt}
%\centering
%  \footnotesize
\caption{Refinement strategies having the best results at the intermediate and   final stage.}
\label{tab:F_summary}
\begin{tabular} {|l?l|l|l|l|l|l|l|}
  \hline
  \multicolumn{8}{|l|}{Intermediate stage} \\
  \hline
  Strategy & Iter & No out & No Coef & Max & Average & Average out & Time  \\
\hlinewd{1pt}
  \multicolumn{8}{|l|}{bi-linear} \\
  \hline
  eFB & 7 & 104 & 62\,291 & 0.258 & 0.028 & 0.216 & \cellcolor{Yellow2}0m18s  \\
  \hline
  bRA tk & 14 & 4\,385 & \cellcolor{Yellow2}36\,540 & 0.324 & 0.043 & 0.218 & 0m32s  \\
  \hline
  \cellcolor{Blue2}eFA tn & 13 & 3\,866 & 40\,472 & 0.465 & 0.042 & 0.222 & 0m27s  \\
\hlinewd{1pt}
  \multicolumn{8}{|l|}{Bi-quadratic} \\
  \hline
  eFB & 7 & 21 & 65\,475 & 0.225 & 0.02 & 0.209 & \cellcolor{Yellow2}0m31s  \\
  \hline
  bRA tk & 15 & 4\,435 & \cellcolor{Yellow2}30\,855 & 0.449 & 0.036 & 0.22 & 1m10s  \\
  \hline
  \cellcolor{Blue2}bRB tk & 7 & 2\,501 & 39\,029 & 0.486 & 0.031 & 0.221 & 0m37s   \\
\hlinewd{1pt}
  \multicolumn{8}{|l|}{bi-cubic} \\
  \hline
  eFB & 7 & 23 & 81\,895 & 0.281 & 0.02 & 0.226 & \cellcolor{Yellow2}0m59s  \\
  \hline
  bRA tk & 24 & 2\,544 & \cellcolor{Yellow2}33\,015 & 0.379 & 0.033 & 0.22 & 3m47s \\
  \hline
  \cellcolor{Blue2}bRB tk & 8 & 3\,674 & 38\,581 & 0.463 & 0.032 & 0.225 & 1m25s  \\
  \hline
\end{tabular} \\

\vspace*{\floatsep}

\begin{tabular} {|l?l|l|l|l|l|l|l|l|}
  \hline
  \multicolumn{9}{|l|}{Final stage} \\
  \hline
  Strategy & Iter & No out & No Coef & Max & Average & Time & Tail & Oscillations \\
\hlinewd{1pt}
  \multicolumn{9}{|l|}{bi-linear} \\
  \hline
  bSB & 8 & 0 & 91\,050 & 0.2 & 0.024 & \cellcolor{Yellow2}0m21s & 1 & No \\
  \hline
  bR+eLA tk & 33 & 0 & \cellcolor{Yellow2}42\,503 & 0.2 & 0.041 & 0m54s & 19 & Medium \\
  \hline
  \cellcolor{Green2}eMcB & 15 & 0 & 49\,212 & 0.2 & 0.035 & 0m27s & 8 & Low \\
  \hline
  \cellcolor{Green2}bRA & 23 & 0 & 47\,430 & 0.2 & 0.037 & 0m40s & 10 & Medium \\
\hlinewd{1pt}
  \multicolumn{9}{|l|}{Bi-quadratic} \\
  \hline
  eFB & 8 & 0 & 65\,523 & 0.2 & 0.02 & \cellcolor{Yellow2}0m33s & 1 & No \\
  \hline
  bSB & 7 & 0 & 111\,992 & 0.2 & 0.015 & \cellcolor{Yellow2}0m33s & 0 & No \\
  \hline
  bR+eLA tk & 31 & 0 & \cellcolor{Yellow2}39\,336 & 0.2 & 0.033 & 1m39s & 18 & Medium \\
  \hline
  \cellcolor{Green2}eMcA & 20 & 0 & 41\,952 & 0.2 & 0.031 & 1m9s & 7 & No \\
\hlinewd{1pt}
  \multicolumn{9}{|l|}{bi-cubic} \\
  \hline
  eFB & 9 & 0 & 82\,064 & 0.199 & 0.02 & \cellcolor{Yellow2}1m6s & 2 & No \\
  \hline
  bRB tk & 40 & 0 & \cellcolor{Yellow2}45\,826 & 0.2 & 0.029 & 3m58s & 32 & Low \\
  \hline
  \cellcolor{Green2}bR+eLB tk & 13 & 0 & 55\,706 & 0.2 & 0.028 & 1m41s & 5 & No \\
  \hline
  \cellcolor{Green2}bR+eLA td & 18 & 0 & 51\,488 & 0.2 & 0.02 & 2m3s & 5 & No \\
  \hline
\end{tabular}
\end{table}
Table~\ref{tab:F_summary} shows the result for some strategies which stand out either
because of low time consumption, few coefficients or an overall good behaviour.
The relation between the performances of the various strategies is not constant throughout the
computation and the first part of the table shows some results at the intermediate stage
of the computation while the second part presents final results.
The intermediate stage is defined as: {\sl At least 99.9\% of the points,
  are closer to the surface than the tolerance and the maximum distance between the point cloud and the
  surface is less than 0.5 meters. This implies that at most  5\,055  
  points is further away from the surface than the tolerance.}

The iteration level where the intermediate stage is reached is shown along
with the number of points outside the tolerance and the number of coefficients in the surface. The
maximum and average distances between the surface and the point cloud are shown next, then follows the average
distance for the points with a larger distance than the tolerance. Distances are given in meters.
At the final stage is a tail entry showing the number of iterations between the intermediate stage
and the final stage added along with an entry called Oscillations.
A high tail entry indicates a slow convergence towards the end of the
computations. The distance to the surface can for some points oscillate around the tolerance value, and the maximum
distance may increase before it decreases again. This behaviour is characterized in the last entry
for the final results and possible values are: No, Low, Medium and High.

Refinement strategies with the lowest time consumption and the leanest surface  are highlighted in yellow. These two success criteria do typically not
coincide and an overall best strategy at the intermediate stage is highlighted in blue. At the
final stage  the overall best strategy is highlighted in green. These selections are subjective choices
weighing time and size.

The full span strategy eFB dominates with respect to low computation time, both at the intermediate and
the final stage. The lowest number of coefficients is reached for 
restricted B-spline based strategies with and without element extension and with threshold based on the number of points outside the tolerance (bRA tk, bRB tk and bR+eLB tk). In most
cases  the least amount of coefficients is reached for strategies with alternating parameter
directions. It is no clear pattern  whether or not alternating parameter directions is preferred
when the execution time and the surface size are judged together. Approximation with bi-cubic
polynomials has the highest computational time for the selected best strategies. 

\begin{figure}
\centering
\includegraphics[width=\textwidth]{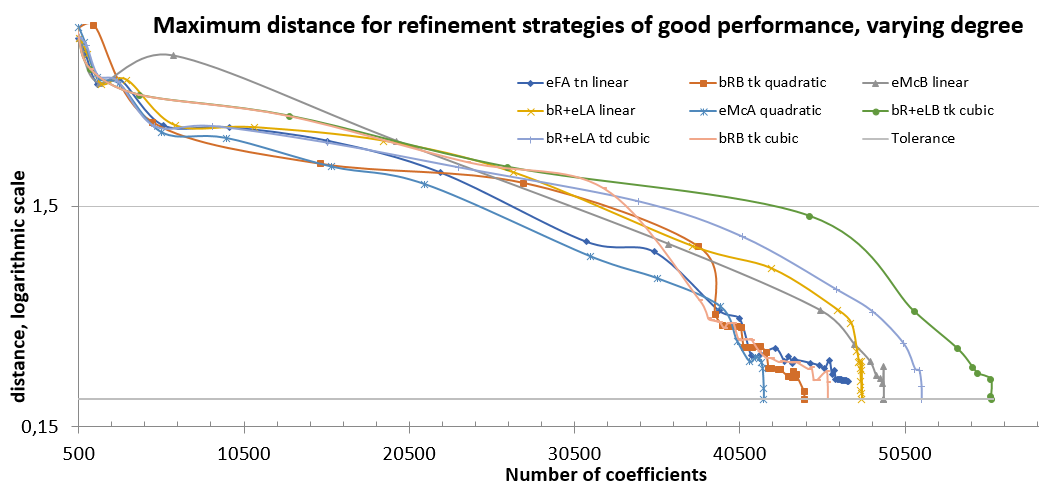}
\caption{Maximum distance for selected well performing refinement strategies, all degrees.}
\label{fig:F_high_dist}

\vspace*{\floatsep}

\includegraphics[width=\textwidth]{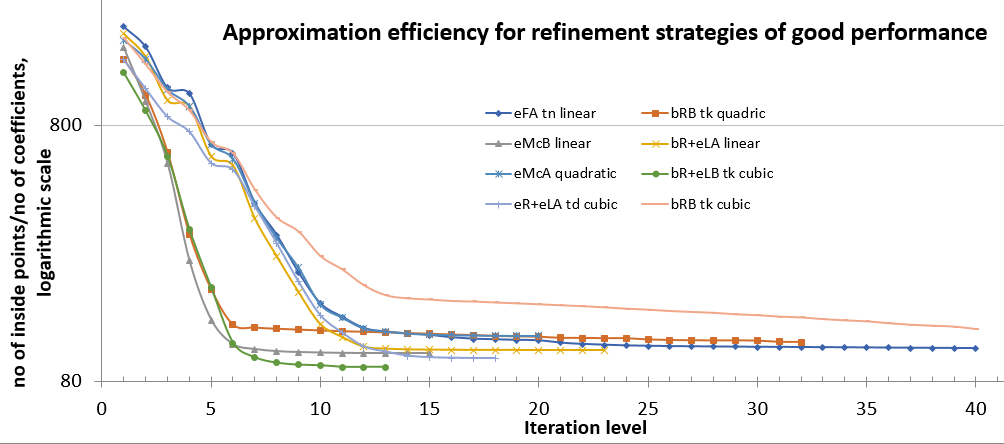}
\caption{Approximation efficiency for for selected well performing refinement strategies, all degrees. A high efficiency
is preferred, but must be coupled with the actual convergence of
the strategy.}
\label{fig:F_high_eff}

\vspace*{\floatsep}

\includegraphics[width=\textwidth]{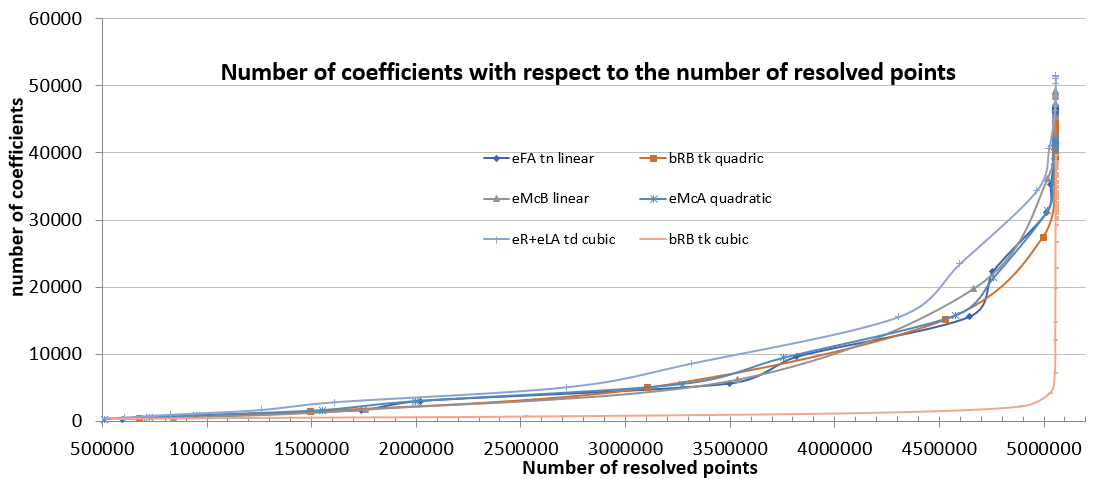}
\caption{Number of coefficients with respect to the number of points
within the tolerance.}
\label{fig:F_high_inside}
\end{figure}
A continuous picture of the performance of the selected best strategies at the intermediate
and final stage is given in
Figures~\ref{fig:F_high_dist} and~\ref{fig:F_high_eff}. Despite a good impression in these figures 
is not bRB tk bi-cubic (restricted B-spline based strategy with a threshold, refines in both parameter
directions) selected as one of the best
strategies for the final result due to the computation time. The strategy has, through parts of the computation, a reasonable maximum
distance between
the surface and the point cloud, it finishes with relatively few coefficients
(Figure~\ref{fig:F_high_dist}) and has a high approximation efficiency (Figure~\ref{fig:F_high_eff}).

The relative performance of the strategies with respect to approximation efficiency is, with the
exception of bRB tk with bi-quadratic polynomials, mainly retained throughout the computations. 
Typically the approximation efficiency is lower for strategies that refine in both parameter directions simultaneously. The increase in the number of coefficients is rapid and the corresponding decrease in the number of unresolved points is not comparatively high.
The ranking of the strategies with respect to the maximum distance varies throughout the iterations. For instance eMcB with bi-linear polynomials (minimum span, combined criterion for selecting the overlapping B-spline,
refines in both parameter directions) increases
the maximum distance early in the computation and continue with a high distance until it drops
below the distance for several other strategies in the last part of the computation.

Minimal span strategy eMcA tn with bi-quadratic polynomials has the best performance throughout the computation. The maximum
distance is constant in the lower part of the group and the approximation efficiency lies in the
group with highest efficiency if bRB tk bi-cubic is kept out. Furthermore, it finishes well within
the allowed number of iterations and with an acceptable computation time. bRB tk with quadratic
polynomials is selected at the intermediate stage and finishes with a good result, but has a relatively poor performance in the first part of the computation. 

Figure~\ref{fig:F_high_inside}
illustrates how the number of coefficients increases with the number
of points inside the tolerance belt for some of the strategies. 
Here, the two bi-cubic cases stand out.
The restricted mesh strategy with element extension has many coefficients
throughout the computation while the version without element extension is very
lean until there is a drastic increase in the number of coefficients when
most points are resolved. The strategies have also different types of thresholds. This is the main source to the difference in the first
part of the computation. In general a
threshold with respect to the number of out-of-tolerance points will lead to less
coefficients with regard to the number of points resolved.

\subsubsection{Other refinement strategies}
Taking also other refinement strategies into account, it is clear that bi-linear approximation in general has lower execution times, but does not necessarily give a
more lean result than bi-quadratic approximation. Bi-cubic approximation leads to higher execution times
and more coefficients. However, the difference between strategies internal to each polynomial degree is in general larger than the
difference between the polynomial degrees.
The selection of degree must also depend on the expected post approximation use of
the surface. Analysis tools like computation of contour curves and evaluation of slope favours a
smooth surface. This test case indicates that the use bi-quadratic LR B-spline surfaces should be
recommended if there is no particular reason for another choice.

\begin{figure}
\centering
\includegraphics[width=\textwidth]{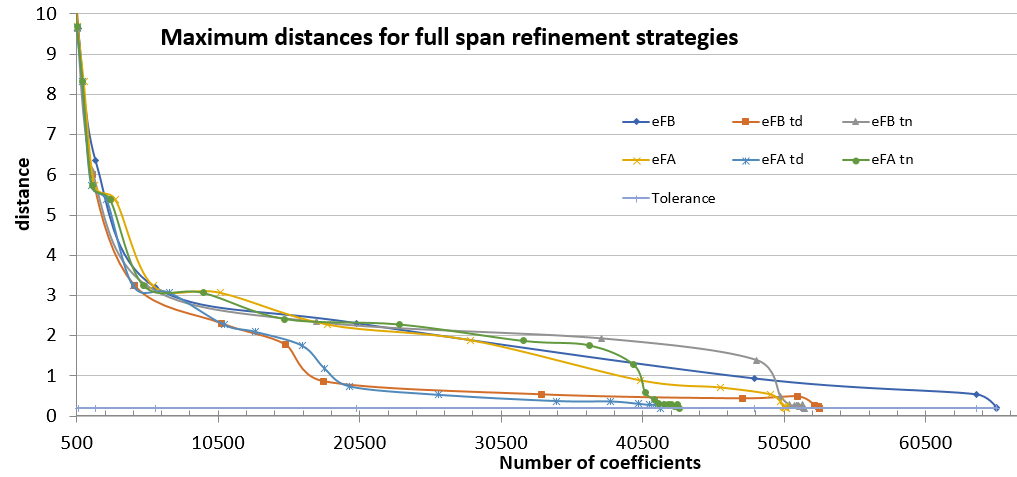}
\caption{Maximum distance for full span strategies with and without threshold, bi-quadratic.}
\label{fig:F_eFB_e2_dist}
\end{figure}
%\vspace*{\floatsep}

\begin{figure}
\centering
\includegraphics[width=\textwidth]{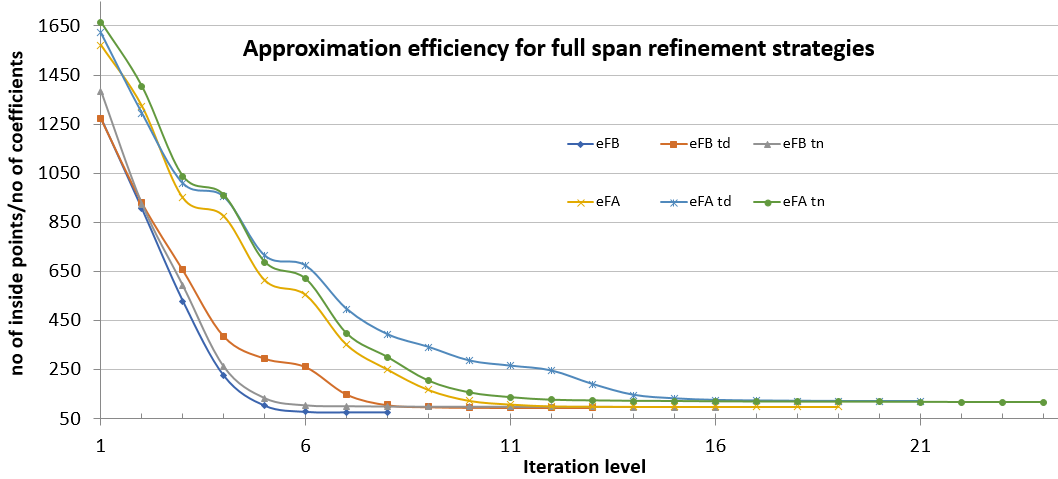}
\caption{Approximation efficiency for full span strategies with and without threshold, bi-quadratic.}
\label{fig:F_eFB_e2_eff}
\end{figure}
Figures~\ref{fig:F_eFB_e2_dist} and~\ref{fig:F_eFB_e2_eff} compare the development of the maximum
distance between the surface and the point cloud, and the approximation efficiency for the full span
refinement strategies with and without a threshold. The maximum distance is seen as a function of the
number of coefficients in the surface, while the approximation efficiency relates to the iteration
level. Bi-quadratic approximation is applied. We see, in both figures, that the graphs of the
refinement strategies gather into two groups, but the groups are not the same. The maximum distance (Figure~\ref{fig:F_eFB_e2_dist})
is reduced more quickly for strategies combined with threshold type td (refine only in elements where the distance between the point cloud 
and the surface is relatively high) than for other combinations.
The evolution for  other strategies diverge at a later stage of the computations.  Strategies with
alternating parameter directions finish with fewer coefficients than strategies that refine
in both parameter direction at every iteration step. Figure~\ref{fig:F_eFB_e2_eff} shows that
strategies with alternating parameter directions have higher approximation efficiency than the others.
Within each group, the strategy that is combined with threshold type td has the highest
efficiency at each iteration level. However, the number of iterations is often higher than for other selections of threshold.
Thus, the td threshold may be beneficial if the computation is stopped at an early stage,
but not necessarily if sufficient iterations to reach convergence is allowed.

Refinement with alternating parameter directions, will in general,  result in less coefficients
than if refinement is applied in both parameter directions at the same iteration level.
This is independent on whether or not a threshold is applied, the type of threshold
and the polynomial degree. This observation is also valid at the intermediate stage
as defined above. Refinement in alternating parameter directions lead to more iterations
although to a lesser extent than we expected based on the fact that the number of
knotline segment inserted at each step is reduced.
The increased number of iterations is only to some extent reflected in higher computation times.

\begin{table}
%\setlength{\tabcolsep}{3.2pt}
%\centering
  %\footnotesize
\caption{Bi-quadratic approximation, structured mesh refinement. Status at the intermediate and final stage. Note that the bSB strategy produces spline spaces similar to those generated by Truncated Hierarchical B-splines (THB). Thus the table gives and indication of the behaviour of the algorithm if THB is used instead of LRB for the bSB strategy. We say "similar" as THB allows narrower refinements and allowing this would probably produce a different spline space  but with similar results.}
\label{tab:Four_bS}
\begin{tabular} {|l?l|l|l|l|l?l|l|l|l|l|}
\hline
Refine & \multicolumn{5}{l?}{Intermediate stage} & \multicolumn{5}{l|}{Final results} \\
strategy & Iter & Pt out & No coef  & Max  & Time &  Iter &
Pt out & No coef  & Max  & Time\\
\hlinewd{1pt}
bSB & 7 & 0 &  \cellcolor{Red2}111\,992 & 0.2 & 0m33s & 7 & 0 & \cellcolor{Red2}111\,992 & 0.2 & \cellcolor{Yellow2}0m33s \\
\hline
bSB td & 7 & 1\,483 & 88\,527 & 0.306 & 0m37s & 10 & 0 & 95\,428 & 0.2 & 0m44s \\
\hlinewd{1pt}
bSA & 13 & 97 & 78\,870 & 0.367 & 0m48s & 15 & 0 & 79\,273 & 0.2 & 0m51s \\
\hline
bSA td & 14 & 646 & 73\,852 & 0.277 & 1m & 18 & 0 & 75\,712 & 0.2 & 1m8s \\
\hline
\end{tabular}
\end{table}

The maximum number of coefficients at the final stage is obtained for the structured mesh
strategy bSB for all
polynomial bi-degrees tested, see Table~\ref{tab:F_range}.  Table~\ref{tab:Four_bS} shows more results for the bi-quadratic case. The numbers highlighted in red tell that this 
strategy obtains the maximum size surface of all strategies. In the bi-quadratic and bi-cubic case, the number of coefficients
for bSB is almost the double of the number for the full span strategy eFB, which also results in large surfaces
compared to other strategies. bSB is expensive also at the
intermediate stage and when run with thresholding td. Similarly, bSA results in more coefficients than any other strategy,
alternating parameter directions or not, except bSB. bSB refines every knot interval in
every B-spline that has any out-of-tolerance points in its domain so many inserted
knotline segments have little impact on the accuracy. The computation times for the
structured mesh strategies are amongst the best, both at the intermediate and at the final
stage.

\begin{table}
%\setlength{\tabcolsep}{3.2pt}
%\centering
  %\footnotesize
\caption{Bi-quadratic approximation, minimum span refinement. Status at the intermediate and final stage.}
\label{tab:Four_eMu}
\begin{tabular} {|l?l|l|l|l|l?l|l|l|l|l|}
\hline
Refine & \multicolumn{5}{l?}{Intermediate stage} & \multicolumn{5}{l|}{Final results} \\
strategy & Iter & Pt out & No coef  & Max  & Time &  Iter &
Pt out & No coef  & Max  & Time\\
\hlinewd{1pt}
eMuB & 7 & 659 & 47\,732 & 0.323& 0m36s & 30 & 3 & 50\,860 & 0.216 & 1.27s \\
\hline
eMuB td & 10 & 3\,482 & 50\,829 & 0.446 & 1m6s & 37 & 6 & 63\,510 & 0.232 & 2m34s \\
\hline
eMuB tn & 11 & 682 & 46\,641 & 0.352 & 0m53s & 40 & 4 & 53\,074 & 0.218 & 2m5s \\
\hlinewd{1pt}
eMuA & 13 & 1\,725 & 39\,933 & 0.4 & 0m55s & 40 & 7 & 45\,762 & 0.224 & 1m57s \\
\hline
eMuA td & 19 & 4\,849 & 44\,017 & 0.412 & 1m41s & 40 & 30 & 57\,771 & 0.235 & \cellcolor{Red2}2m58s \\
\hline
eMuA tn & 14 & 4\,069 & 37\,565 & 0.392 & 1m5s & 40 & 62 & 54\,101 & 0.239 & 2m29s \\
\hline
\end{tabular}
\end{table}

The minimum span strategies with sub type u (select the overlapping B-spline with the relative highest number out-of-tolerance points in its domain for refinement) fail to converge for all  bi-degrees, with and
without threshold. The figures for the bi-quadratic case are shown
in Table~\ref{tab:Four_eMu}. The same strategies do fairly well at the intermediate stage and also
for some more iterations. The number of coefficients are kept low, the computation times
are higher than most other methods, but not with a significant amount at that
stage. eMuA/eMuB is a minimum span strategy with a quite restrictive selection of elements to
refine. This implies that only the most important knot insertions are performed at each
iteration step. There is a risk of blocking further convergence by failing to split
an important element. Blocking does indeed occur, even for this relatively
smooth data set. The minimum span strategies eMlA and eMlB, which select to refine B-splines with the maximum support, have an average performance. eMcB is a combination of eMlB and eMuB, and eMcA is
the corresponding strategy with alternating parameter directions. These strategies show
good results.

The restricted B-spline based strategies bRB and bRA with associated threshold tend to
give good results at the intermediate stage, but have not consistently converged at the
final stage. The refinement restrictions become exaggerated, in particular when
threshold tk or td+k, which combines td and tk, are applied. Adding the element based
fall-back embedded in the strategies bR+eLA and bR+eLB give lower execution times, consistent
convergence and the surface size is kept at the same level or reduced.

\subsection{Gaustatoppen}
\subsubsection{Data set}
\begin{figure}
  \centering
\begin{tabular}{cc}
\includegraphics[width=0.46\textwidth]{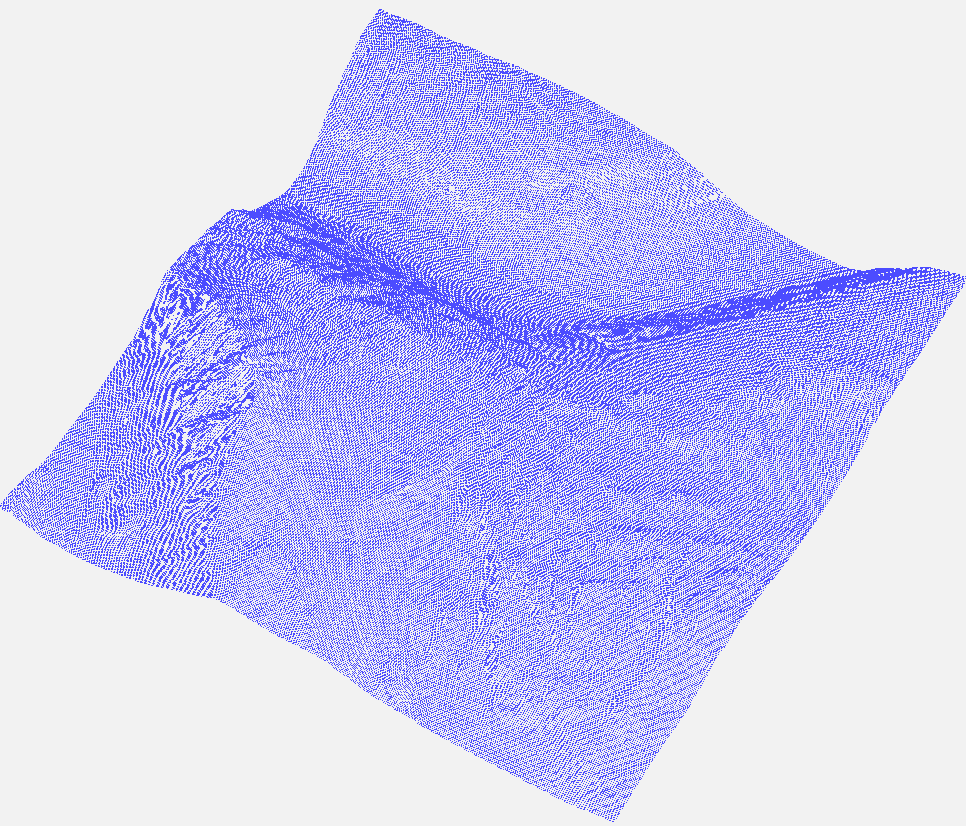}
&\includegraphics[width=0.43\textwidth]{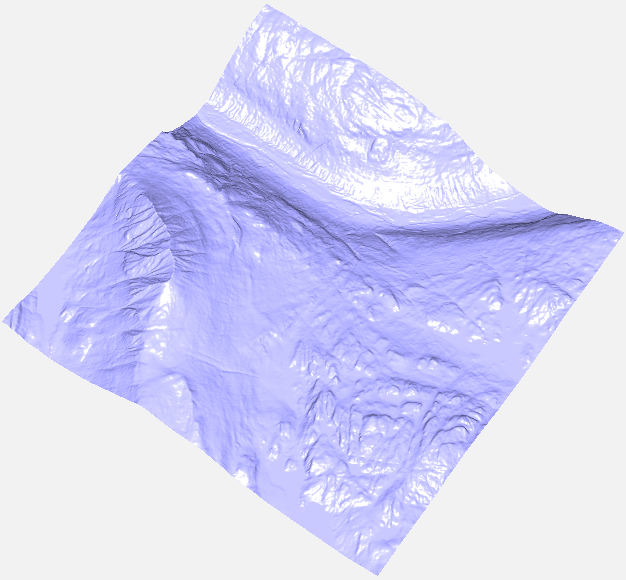} \\
(a) & (b)
\end{tabular}
\caption{(a) Raster points from Gaustatoppen, Norway. (b) Corresponding LR B-spline
  surface. \label{fig:gaustapts}}
\end{figure}
Figure~\ref{fig:gaustapts} (a) shows a completely regular point set containing 490\,000
points. The points are extracted from a sparse raster approximating measurement data
from the mountain area Gaustatoppen in Norway. The elevation range is
$[216.6, 1877.3]$ measured in meters and the standard deviation is 370.542.
The data set covers an area of 6\,507 times 6\,507 square kilometers. The data set is regular, but
very sparse. It has already been
processed and thinned considerably. The distances between neighbouring points are 13 meters and 
being extracted from a raster, the point cloud represents a piecewise constant surface.
As every point carries a lot of information, an accurate approximation of the points
with a lean surface is not feasible. However, the data set can still be used to
distinguish between different refinement strategies, and in particular evaluate the strategies
with respect to robustness.

\subsubsection{Result overview}
The Gaustatoppen terrain is much more demanding than in the previous case. The
area has one high mountain top and a valley with a river. There are also several other tops and
valleys in the neighbourhood. 
We specify a tolerance of 0.5 meters and a maximum number of iterations of 40.
The intermediate stage is defined as the iteration level where
{\sl the maximum distance between the point cloud and the
surface is less than two meters, the number of points where the tolerance is not satisfied does not
exceed 4\,900 and the average distance in out-of-tolerance points is less than one meter}.
\begin{table}
  \caption{Computation time and number of coefficients. Performance
  ranges and associated strategies for polynomial degrees one to three. }
  \label{tab:G_range}
  \begin{tabular}{|l||l|l||l|l||l|l||l|l|}
    \hline
Deg & Min t. & Strategy & Max t. & Strategy & Min cf & Strategy & Max cf & Strategy \\
\hline
\hline
1 & 1m26s & eFA & 29m402 & eMuB td & 270\,613 & bR+eLA tk & 1\,227\,085 & bSB \\
\hline
2 & 0m31s & eFA & 1h23m28s & eMuB td & 195\,532 & bRA tk & 814\,295 & bSB \\
\hline
3 & 0m52s & bSB & 3h13m30s & eMuB td & 233\,482 & bRA td+k & 757\,012 & eMuB td \\
\hline
\end{tabular}
\end{table}

There exists a turning point where the surface size and the distribution of new knotlines cause the
time spent in maintaining the data structure of the surface during refining to increase very rapidly.
This point is reached for some refinement
strategies for this test case. Table~\ref{tab:G_range} presents the
range in computational time and final surface size for all polygonal
degrees. 
The highest execution
times occur for minimum span strategies with
subtype u and threshold based on distance (td).
The minimum execution times are obtained for the
full span strategy with alternating parameter directions and the
structured mesh strategy with refinement in both directions. The latter
strategy occurs also as the one that results in most coefficients and
contrary to the normal pattern do bi-linear polynomials lead to
the largest surface size. The least surface size is obtained for various
restricted mesh strategies with threshold.

Several refinement strategies fail to resolve all points and
also to reach the intermediate stage. The required maximum distance of less than 2 meters 
is most demanding. Bi-linear approximation struggles if a threshold is applied. A restraint
selection of elements to refine leads to convergence
failure also without additional thresholds. In particular, most
minimum span element based refinement strategies fail to reach the required
accuracy within the specified number of iterations. The full span strategies converge in
most cases, the exception is eFA td (alternating parameter directions and
threshold with respect to distance). B-spline based strategies converge more frequently, but also
these strategies are unstable, especially when a threshold is applied.

The bi-quadratic and bi-cubic approximation have similar convergence patterns as in the bi-linear case, but less extreme. More strategies reach the intermediate stage.
The bi-cubic case results in larger surfaces and has higher execution times than the bi-quadratic case.

\subsubsection{Best performances}
\begin{table}
%\setlength{\tabcolsep}{4pt}
%\centering
%  \footnotesize
\caption{Refinement strategies having the best results at the intermediate and
  the final stage.}
\label{tab:G_summary}
\begin{tabular} {|l?l|l|l|l|l|l|l|}
  \hline
  \multicolumn{8}{|l|}{Intermediate stage} \\
  \hline
  Strategy & Iter & No out & No Coef & Max & Average & Average out & Time  \\
\hlinewd{1pt}
  \multicolumn{8}{|l|}{bi-linear} \\
  \hline
  bR+eLA & 30 & 15 & 265\,776 & 1.73 & 0.113 & 0.971 & \cellcolor{Yellow2}1m13s  \\
  \hline
  \cellcolor{Blue2}bRA tk & 39 & 1 & \cellcolor{Yellow2}265\,050 & 1.788 & 0.118 & 1.788 & 1m47s \\
\hlinewd{1pt}
  \multicolumn{8}{|l|}{Bi-quadratic} \\
  \hline
  eFB tn & 6 & 2\,942 & 214\,016 & 1.452 & 0.106 & 0.591 & \cellcolor{Yellow2}0m13s  \\
  \hline
  \cellcolor{Blue2}bRA tk & 13 & 4\,006 & \cellcolor{Yellow2}171\,341 & 1.681 & 0.13 & 0.588 & 0m24s  \\
\hlinewd{1pt}
  \multicolumn{8}{|l|}{bi-cubic} \\
  \hline
  eMlB tn & 7 & 1\,953 & 261\,332 & 1.762 & 0.097 & 0.583 & \cellcolor{Yellow2}0m21s  \\
  \hline
  \cellcolor{Blue2}bRB tk & 7 & 4\,797 & \cellcolor{Yellow2}193\,719 & 1.908 & 0.121 & 0.602 & 0m30s  \\
  \hline
\end{tabular}

\vspace*{\floatsep}

\begin{tabular} {|l?l|l|l|l|l|l|l|l|}
  \hline
  \multicolumn{9}{|l|}{Final stage} \\
  \hline
  Strategy & Iter & No out & No Coef & Max & Average & Time & Tail & Oscillations \\
\hlinewd{1pt}
  \multicolumn{9}{|l|}{bi-linear} \\
  \hline
  \cellcolor{Green2}eFA & 38 & 0 & 271\,450 & 0.5 & 0.112 & \cellcolor{Yellow2}1m26s & 0 & Medium \\
  \hline
  bR+eLA tk & 39 & 0 & \cellcolor{Yellow2}270\,613 & 0.5 & 0.118 & 1m29s & 8 & Medium \\
\hlinewd{1pt}
  \multicolumn{9}{|l|}{Bi-quadratic} \\
  \hline
  eFA & 18 & 0 & 213\,774 & 0.5 & 0.099 & \cellcolor{Yellow2}0m31s & 6 & Low \\
  \hline
  \cellcolor{Green2}eFA tn & 22 & 0 & 198\,013 & 0.5 & 0.11 & 0m41s & 9 & Low \\
  \hline
  \cellcolor{Green2}bRA tk & 22 & 0 & \cellcolor{Yellow2}195\,532 & 0.5 & 0.112 & 0m53s & 9 & Low \\
  \hline
  \cellcolor{Green2}bR+eLA td+k & 19 & 0 & 199\,654 & 0.5 & 0.103 & 0m37s & 7 & Low \\
\hlinewd{1pt}
  \multicolumn{9}{|l|}{bi-cubic} \\
  \hline
  bRB & 9 & 0 & 355\,993 & 0.499 & 0.07 & \cellcolor{Yellow2}0m52s & 2 & No \\
  \hline
  bRA td+k & 35 & 0 & \cellcolor{Yellow2}233\,482 & 0.5 & 0.103 & 2m51s & 13 & Medium \\
  \hline
  \cellcolor{Green2}eFA tn & 21 & 0 & 233\,552 & 0.5 & 0.096 & 1m7s & 9 & Low \\
  \hline
\end{tabular}
\end{table}
Table~\ref{tab:G_summary} presents the best results for the Gaustatoppen test case. bRA tk (restricted mesh, alternating parameter directions,
threshold with respect to number of out-of-tolerance points) with
bi-linear polynomials is selected
as the best result at the intermediate stage despite the large maximum distance since the number of
out-of-tolerance points is only one.  
Independent of the bi-degrees is the strategy with the least number of coefficients at this stage
the restricted B-spline based method with threshold tk, and it is also
selected as the overall best strategy at the intermediate stage.

Most best performing
strategies at the final stage refine in alternating parameter
directions, and strategies with threshold tend to get the best overall score. 

The oscillation entry covers oscillations in both the maximum distance and the number of unresolved
points. Typical for the selected strategies at the final
stage is that they have some oscillations in the maximum
distance and hardly any oscillations in the number of outside points. The selected strategies in the bi-quadratic case outperform the selections in the bi-linear and bi-cubic cases both with regard to execution time and surface size.

\begin{figure}
\centering
\includegraphics[width=\textwidth]{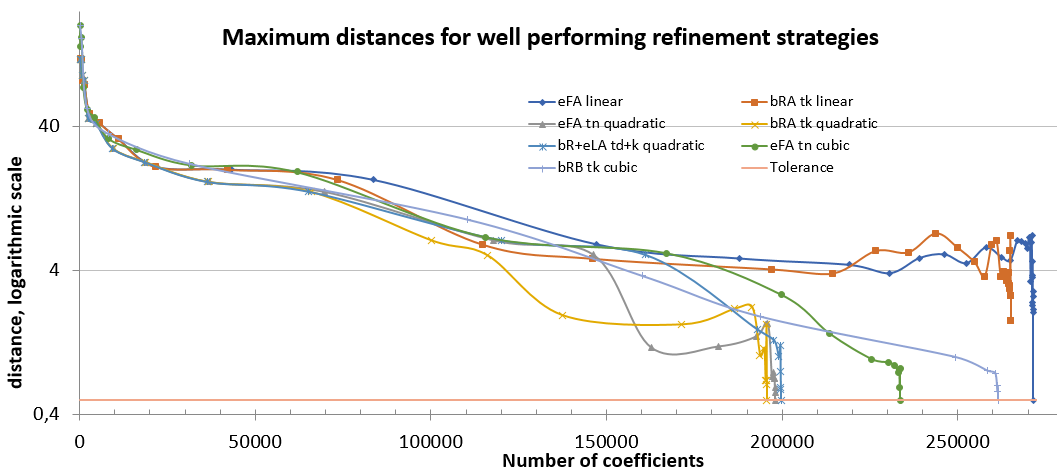}
\caption{Maximum distance for selected well performing refinement strategies, all degrees.}
\label{fig:G_high_dist}
\end{figure}
%\vspace*{\floatsep}

\begin{figure}
\centering
\includegraphics[width=\textwidth]{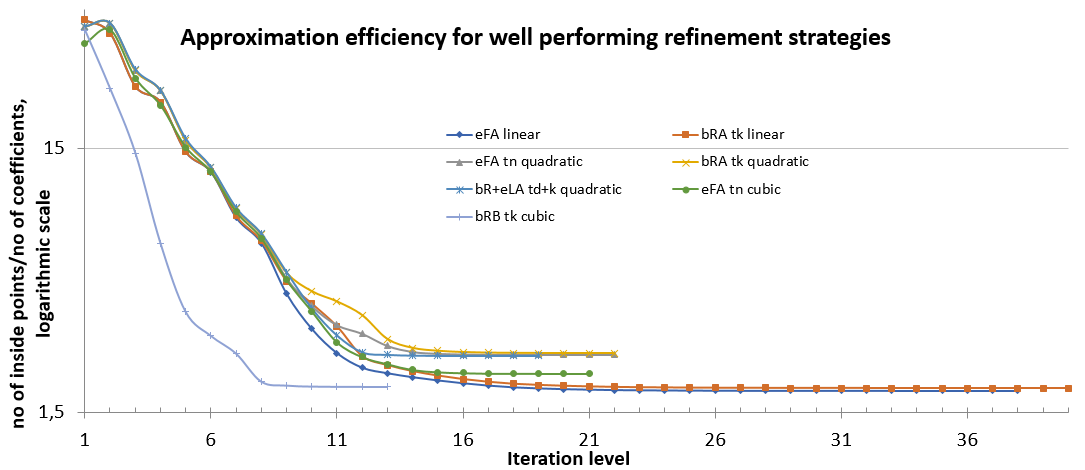}
\caption{Approximation efficiency for selected well performing refinement strategies, all degrees.}
\label{fig:G_high_eff}
\end{figure}
B-spline based strategies dominate at the intermediate stage, while at the
final stage are element based strategies better when surface
size and execution time are viewed together.
Figures~\ref{fig:G_high_dist} and~\ref{fig:G_high_eff} provide more information on the behaviour of
the selected best strategies throughout the entire computations. The maximum distance decreases steadily
for a while before it for most strategies starts to oscillate. The oscillations are strongest in the bi-linear case and 
least profound for bi-cubic polynomials. The maximum distance is, over most of the computations, lowest
with respect to the number of coefficients for bRA tk bi-quadratic (restricted mesh, threshold with respect to out-of-tolerance points), This
combination also has the highest approximation efficiency, see
Figure~\ref{fig:G_high_eff}. The strategy that refines in both
parameter directions simultaneously (bRB tk bi-cubic) has a
steeper descent in approximation efficiency, but ends up at the
same level as several of the other strategies and the resulting
number of coefficients (261\,462) is quite competitive.

\subsubsection{Other refinement strategies}
Approximations with bi-quadratic or bi-cubic polynomials are more robust than  bi-linear approximation.
However, the minimum span refinement
strategies struggles to converge in many cases also for degree two and three. eMlB, which refines the largest B-spline overlapping
a selected elements, converges both with and without threshold, but the execution time for eMlB td
is high. eMu (the overlapping B-spline with most unresolved points is split) does not converge for any combination of threshold types and polynomial degrees
and fail to reach the intermediate stage for some combinations. The only eMc (combination of size and unresolved points) strategy that converge completely is bi-cubic eMcB without a threshold. The numbers of outside points is, however, low at the finishing stage. The execution time is high for threshold type td, but roughly in
the same magnitude as the full span strategies for other cases. 

\begin{table}
%\setlength{\tabcolsep}{3.2pt}
%\centering
  %\footnotesize
\caption{Bi-quadratic approximation, combined refinement strategies. Status at the intermediate and final stage.}
\label{tab:Gaus_combined}
\begin{tabular} {|l?l|l|l|l|l?l|l|l|l|l|}
\hline
Refine & \multicolumn{5}{l?}{Intermediate stage} & \multicolumn{5}{l|}{Final results} \\
strategy & Iter & Pt out & No coef  & Max  & Time &  Iter &
Pt out & No coef  & Max  & Time\\
\hlinewd{1pt}
eMcA tn & 15 & 2\,101 & 194\,789 & 1.642 & 0m43s & 40 & 3 & 211\,543 & 2.527 & 2m20s \\
\hline
eMc/FA tn & 15 & 2\,101 & 194\,789 & 1.642 & 0m43s & 29 & 0 & 211\,336 & 0.5 & 1m40s \\
\hlinewd{1pt}
bRA tk & 13 & 4\,006 & 171\,341 & 1.681 & 0m24s & 22 & 0 & 195\,532 & 0.5 & 0m53s \\
\hline
bR+eLA tk & 12 & 2\,018 & 190\,310 & 1.349 & 0m23s & 20 & 0 & 201\,934 & 0.5 & 0m45s \\
\hline
\end{tabular}
\end{table}

The restricted B-spline based strategies bRA and bRB and their element extended counterparts bR+eLA and
bR+eLB converge for all threshold combinations and polynomial degrees two and three. In the
bi-linear case do most combinations involving threshold with respect to the
number of out-of-tolerance points struggle, but in general do the restricted B-spline based strategies
perform better than the element based ones. 

\begin{figure}
  \begin{center}
\includegraphics[width=\textwidth]{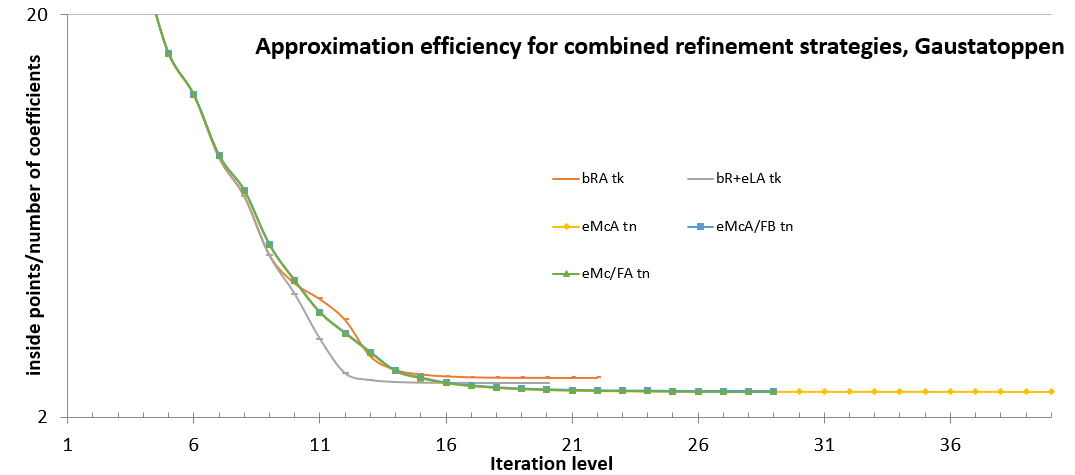}
  \caption{Approximation efficiency for selected combined refinement strategies with threshold type tk or tn, bi-quadratic polynomials.
    \label{fig:Gausta_combined}}
\end{center}
\end{figure}
We focus now on the
minimum span strategy eMcA with and without a full span back-up and the restricted B-spline based strategy bRA with and without element extension.
Threshold types tn and tk (point based) are applied. Table~\ref{tab:Gaus_combined} gives the numbers at the intermediate and final stage in the bi-quadratic case, and Figure~\ref{fig:Gausta_combined} compares the
approximation efficiency for the middle and last part of the computation.
We see that bRA tk has the best efficiency. bR+eLA tk
finishes second best, but has a lower efficiency in some part of the iterations. eMcA tn does not
converge completely alone, but do so when combined with a
full span strategy. The B-spline based strategies covered here
have better execution times than the minimum span strategies.

The surface size of all structured mesh strategies (bSB, bSB td, bSA and bSA td) is high and
exceeds the size of most other
strategies. For degree two and three are the highest number of coefficients obtained for strategies that totally
fail to converge with minimum span strategy eMuA td 
as the most extreme. The structured mesh strategies and the minimum span strategies
with sub type u turn out to be inappropriate for approximation of large sets of geospatial data points for
opposite reasons. One approach produces surfaces with too much data and the other is too restrictive.
Both groups of strategies will be omitted for the remaining test cases.

For this data set the number of coefficients exceeds the number of data points for several
combinations of refinement strategies and threshold. In general  the minimum element size in the
final surface is smaller than the point cloud resolution. Restricting the element size to this resolution,
i.e. the element cannot exceed 13 meters in any direction, leads to an almost regular spline surface.
Most of the surface is refined to account for the details in the terrain, but the refinement in the
steepest areas is not sufficient to obtain the requested accuracy. This tendency is strongest in the
bi-linear case, but occurs also for bi-quadratic and bi-cubic polynomials. For instance, applying
bi-quadratic  approximation
and the restricted mesh strategy bRA tk with restrictions on the element size result
in a surface with 124\,141 coefficients. The maximum distance between the point cloud and the
surface is 5.629 meters and 26\,844 points have a larger distance than 0.5 meters. The same setup
without element restriction converges with 195\,532 coefficients.

\subsection{Scan lines} \label{sec:scanline}
\subsubsection{Data set}
\begin{figure}
\centering
\begin{tabular}{ccc}
\includegraphics[width=0.35\textwidth]{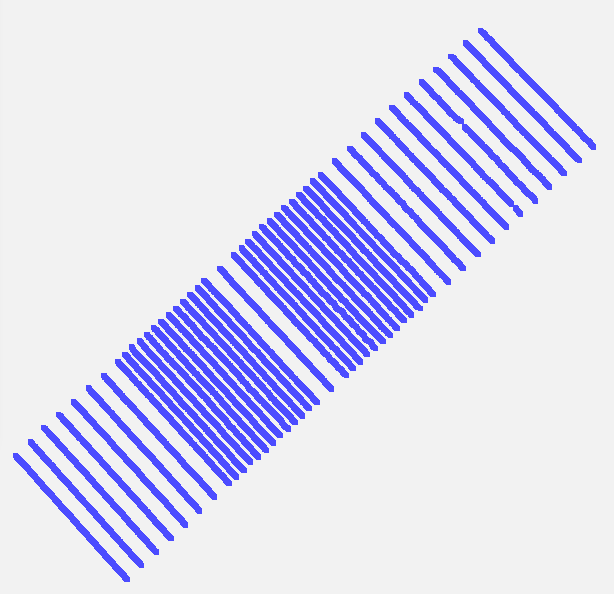} 
&\includegraphics[width=0.27\textwidth]{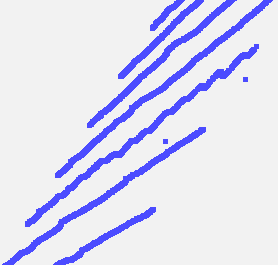} 
&\includegraphics[width=0.35\textwidth]{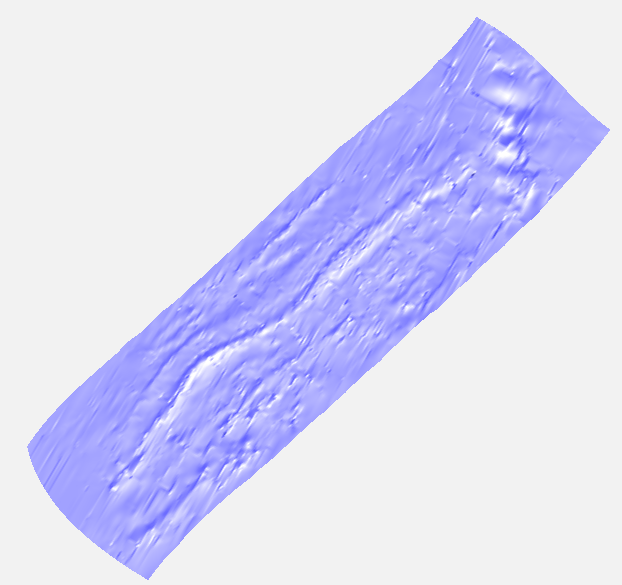} \\
(a) & (b) & (c)
\end{tabular}
\caption{ \label{fig:scanline_pts}
  (a) Sub sea scan line points, English channel. (b) Detail emphasizeing outliers. (c) Corresponding surface.}
\end{figure}
Figure~\ref{fig:scanline_pts} (a) shows a small data set of 71\,888 points. The
data set is obtained from the English channel and the
points are organized in scan lines with varying distances. A detail can be seen in (b). The point set
contains a number of outlier points situated off the scan lines, which are
inconsistent with the general trend in the data. The corresponding surface (c)
is created with an accuracy of 0.5 meters and is approximating also the
outlier points closely.
The regular points represent a smooth
part of the sea bed. The data set covers about 23
square kilometers, the elevation range is
$[-47.8, -12.1]$ meters and the standard deviation is 2.408. 

\subsubsection{Results overview}
The approximation is run with a tolerance of 0.5 meters and the maximum number of iterations is 30
and 40 for refinement strategies refining in both parameter directions at every iteration level and with alternating parameter directions, respectively.
An intermediate stage is defined when 99\% of the points are within the resolution. This
corresponds to 719 unresolved points. Due to the outliers is no condition put on the maximum distance at the
intermediate stage.

\begin{table}
  \caption{Computation time and number of coefficients. Performance
  ranges and associated strategies for polynomial degrees one to three. }
  \label{tab:Sc_range}
  \begin{tabular}{|l||l|l||l|l||l|l||l|l|}
    \hline
Deg & Min t. & Strategy & Max t. & Strategy & Min cf & Strategy & Max cf & Strategy \\
\hline
\hline
1 & 0m1.3s & eFB & 0m4.6s & eMcA tn & 10\,593 & eFA td & 11\,978 & bRA \\
\hline
2 & 0m3.1s & eFB & 0m12.4s & eMlA tn & 16\,453 & eFA td & 23\,264 & bR+eLB tk \\
\hline
3 & 0m8.9s & eFB & 0m35.6s & eMcA & 25\,331 & eFA tn & 35\,386 & bRB td+k \\
\hline
\end{tabular}
\end{table}

The current
test case contains few points and this is reflected in fast execution times also for the
strategies with the poorest performance as can be seen in Table~\ref{tab:Sc_range}. The full span strategy with refinement in
both parameter direction is fastest while minimum span strategies
require most time and do not always converge within the given
number of iterations. Most strategies do converge, some would converge if more iterations
were allowed and a few fail completely. A high number of iterations tend to correspond with high execution times. 
Full span strategies with alternating parameter
directions result in the least number of coefficients and 
restricted mesh strategies in the highest.
Bi-linear approximation
results in the smallest number of coefficients and the least execution times while
bi-cubic gives higher execution times and more coefficients. The number of strategies that
do not converge during the specified number of iterations is the least in the quadratic
case.

\subsubsection{Best performances}
\begin{table}
%\setlength{\tabcolsep}{4pt}
%\centering
%\footnotesize
\caption{Scan lines, refinement strategies with the best results.}
\label{tab:Sc_summary}
\begin{tabular} {|l?l|l|l|l|l|l|l|}
  \hline
  \multicolumn{8}{|l|}{Intermediate stage} \\
  \hline
  Strategy & It & No out & No Coef & Max & Average & Average out & Time  \\
\hlinewd{1pt}
  \multicolumn{8}{|l|}{bi-linear} \\
  \hline
  eFB & 5 & 587 & 11\,001 & 25.235 & 0.117 & 1.011 & \cellcolor{Yellow2}0m0.5s  \\
  \hline
  eFB tn & 6 & 604 & 7\,574 & 17.377 & 0.124 & 0.921 & \cellcolor{Yellow2}0m0.5s  \\
  \hline
  eMcB & 5 & 709 & 8\,792 & 25\,249 & 0.123 & 0.963 & \cellcolor{Yellow2}0m0.5s  \\
  \hline
  eMcB tn & 6 & 658 & 6\,483 & 23.884 & 0.126 & 0.933 & \cellcolor{Yellow2}0m0.5s  \\
  \hline
  bRB & 5 & 612 & 10\,853 & 25.235 & 0.119 & 0.995 & \cellcolor{Yellow2}0m0.5s  \\
  \hline
  \cellcolor{Blue2}eMlB & 6 & 617 & 6\,884 & 18.22 & 0.126 & 0.949 & \cellcolor{Yellow2}0m0.5s  \\
  \hline
  \cellcolor{Blue2}eMcA tn & 13 & 672 & \cellcolor{Yellow2}6\,313 & 19.893 & 0.129 & 0.885 & 0m0.9s  \\
\hlinewd{1pt}
  \multicolumn{8}{|l|}{Bi-quadratic} \\
  \hline
  eFB & 5 & 606 & 14\,422 & 23.663 & 0.114 & 1.071 & \cellcolor{Yellow2}0m1s  \\
  \hline
  eFB tn & 6 & 603 & 10\,074 & 18.220 & 0.118 & 1.015  & \cellcolor{Yellow2}0m1s \\
  \hline
  \cellcolor{Blue2}bRB td+k & 22 & 667 & \cellcolor{Yellow2}8\,235 & 13.404 & 0.124 & 0.699 & 0m2.5s  \\
  \hline
  \cellcolor{Blue2}eFA tn & 12 & 664 & 8\,848 & 18.153 & 0.123 & 0.962 & 0m5.3s  \\
\hlinewd{1pt}
  \multicolumn{8}{|l|}{bi-cubic} \\
  \hline
  eFB tn & 6 & 614 & 12\,387 & 18.919 & 0.117 & 1.023 & \cellcolor{Yellow2}0m1.8s  \\
  \hline
  bRA td+k & 34 & 696 & \cellcolor{Yellow2}9\,955 & 17.748 & 0.125 & 0.851 & 0m7s  \\
  \hline
  \cellcolor{Blue2}eFA tn & 12 & 635 & 11\,231 & 18.964 & 0.12 & 1.0 & 0m2.5s  \\
  \hline
  \cellcolor{Blue2}bRA tk & 22 & 712 & 9\,994 & 23.505 & 0.126 & 0.938 & 0m5.1s  \\
  \hline
  \cellcolor{Blue2}bR/eFA tn/k & 22 & 712 & 9\,994 & 23.505 & 0.126 & 0.938 & 0m5.1s  \\
\hline
\end{tabular} 

\vspace*{\floatsep}

\begin{tabular} {|l?l|l|l|l|l|l|l|l|}
  \hline
  \multicolumn{9}{|l|}{Final stage} \\
  \hline
  Strategy & It & No out & No Coef & Max & Average & Time & Tail & Oscillations \\
\hlinewd{1pt}
  \multicolumn{9}{|l|}{bi-linear} \\
  \hline
  eFB & 10 & 0 & 14\,678 & 0.5 & 0.108 & \cellcolor{Yellow2}0m1.3s & 5 & No \\
  \hline
  eFA td & 36 & 0 & \cellcolor{Yellow2}10\,593 & 0.5 & 0.113 & 0m2.4s & 12 & Medium \\
  \hline
  \cellcolor{Green2}bRA & 18 & 0 & 11\,978 & 0.5 & 0.113 & 0m1.6s & 8 & No \\
\hlinewd{1pt}
  \multicolumn{9}{|l|}{Bi-quadratic} \\
  \hline
  eFB & 11 & 0 & 21\,058 & 0.5 & 0.102 & \cellcolor{Yellow2}0m3.1s & 6 & No \\
  \hline
  \cellcolor{Green2}eFA td & 35 & 0 & \cellcolor{Yellow2}16\,453 & 0.5 & 0.105 & 0m5.6s & 11 & Medium \\ 
\hlinewd{1pt}
  \multicolumn{9}{|l|}{bi-cubic} \\
  \hline
  eFB & 13 & 0 & 29\,908 & 0.499 & 0.099 & \cellcolor{Yellow2}0m8.9s & 7 & Low \\
  \hline
  \cellcolor{Green2}eFA tn & 24 & 0 & \cellcolor{Yellow2}25\,331 & 0.5 & 0.103 & 0m10.4s & 12 & Low \\
  \hline
\end{tabular} 
\end{table}
Table~\ref{tab:Sc_summary} shows the refinement strategies that perform best with regard
to execution time, number of coefficients and  overall result.
The first part presents the best strategies at the intermediate stage. Entries marked with
yellow indicate the best result for one single category and we see that
several strategies have the same execution time in the bi-linear case. 
The overall best strategies are  minimum
span with different sub categories, one refines in both parameter directions at each iteration step and the other in one direction. 

Various full span strategies and one restricted mesh strategy stand out when bi-quadratic polynomials are applied. In the bi-cubic case, three strategies are selected.
bRA tk and bR/eFA tn/k are identical at the intermediate stage. The first is a restricted mesh strategy and the second a combination of
a restricted mesh and a full span strategy, and the switch to eFA has not yet taken place. The full span strategy eFA tn has more
coefficients than the other selected bi-cubic strategies, but has a lower execution time.

The second part of the table presents the best results at the final stage. The full span strategy with refinement in both
parameter directions at once has the fastest execution time
for all polynomial degrees. The full span strategy with alternating
parameter directions (eFA) and some threshold results in the smallest
surfaces and is also selected as the overall best strategy in
the quadratic and cubic case. In the linear case is a
restricted mesh strategy selected.
eFA tn bi-linear and eFA td bi-quadratic are given the oscillation flag medium due to
the behaviour of the maximum
distance. The number of outside points does not oscillate although it decreases slowly in the
last few iterations.
Strategies with alternating and not alternating
parameter directions are chosen as best at the
intermediate stage, but at the final stage do 
strategies with refinement in both directions occur only as the ones with lowest computation time.

\begin{figure}
\centering
\includegraphics[width=0.9\textwidth]{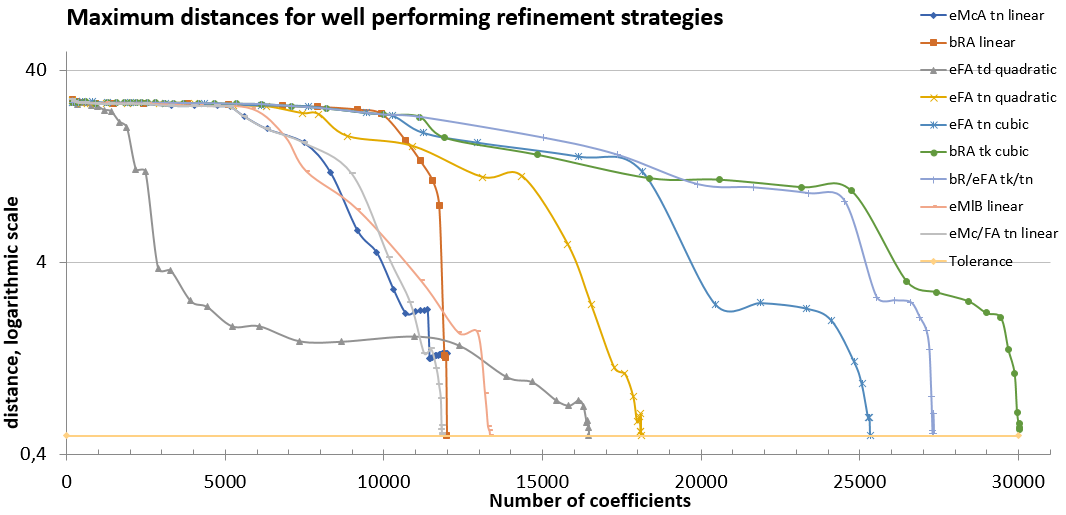}
\caption{The scanline data set.
  Maximum distance for selected well performing refinement strategies, all degrees.}
\label{fig:Sc_high_dist}
\end{figure}
%\vspace*{\floatsep}
\begin{figure}
\includegraphics[width=\textwidth]{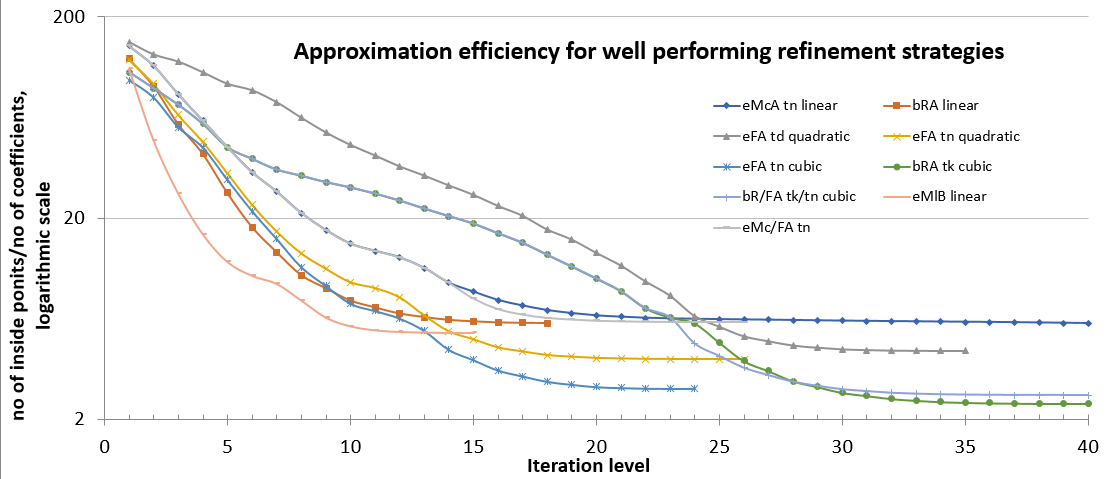}
\caption{Approximation efficiency for selected well performing refinement strategies, all degrees.}
\label{fig:Sc_high_eff}
\end{figure}
The strategies that perform best at the intermediate stage are not always those that finish with
the best result. The Figures~\ref{fig:Sc_high_dist} and~\ref{fig:Sc_high_eff} give a more
continuous picture of the behaviour of the selected strategies. Approximation with bi-cubic
polynomials results in most coefficients and the mainly B-spline based bi-cubic strategies
have a low approximation efficiency in the later part of the iterations. Early in the
iterations do these two strategies have a good efficiency. Bi-linear approximation gives the
least number of coefficients. The minimum span strategy eMcA tn bi-linear does not converge in isolation, but 
changing to eFA tn when the convergence speed decreases results in convergence and a moderate
number of coefficients. Most strategies maintain a high maximum distance in the start
before it starts decreasing and in particular the restricted mesh strategy bRA bi-linear has a very rapid decrease. This
is the method that requires the least number of iterations.

The strategy with a 
threshold based on distance stands out (eFA td quadratic). The maximum distance decreases fast and the approximation
efficiency is high for most of the iterations. This pattern can be found when this type of threshold is
applied also for several other refinement strategies and polynomial degrees. On the other hand, the
combination of element based strategies and td threshold frequently leads to lack of convergence
or slow convergence in the last part of the computations, and the number of points outside the
tolerance belt decreases more rapidly without threshold or if a threshold based on the number of out-of-tolerance points (tn or tk) is
applied. Type td threshold has a less aggressive introduction of coefficients than other methods.
This influences the number of unresolved points as well as the approximation efficiency.

The difference between the restricted B-spline based strategies bRB/bRA and corresponding element extended
strategies is moderate with respect to data size and execution time, but the pure restriced mesh strategies with
threshold have occurrences of slow or missing convergence. Neither of the strategies compete with
the full span element based method with regard to data size and performance although bRA tk is selected
as one of the preferred strategies at the intermediate stage for bi-cubic polynomials.

The scan line data set is, despite the outlier points,  approximated with good accuracy at the cost of
surface smoothness. An alternative stop criteria
for configurations where particular points have a severe lack of convergence should be considered.

\subsection{Sore Sunnm\o re, sea bed} \label{sec:Sore}
\subsubsection{Data set} 
\begin{figure}
  %\centering
\begin{tabular}{cc}
\includegraphics[width=0.48\textwidth]{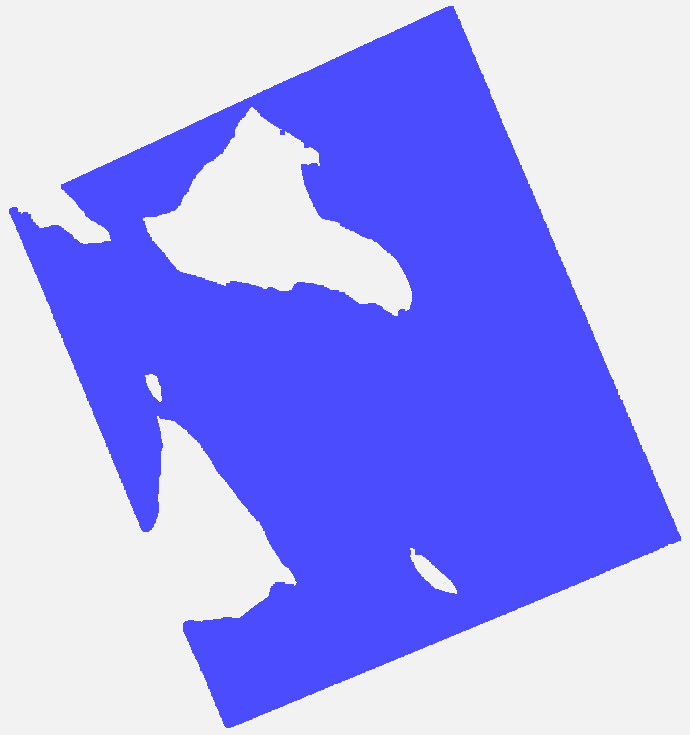}
&\includegraphics[width=0.48\textwidth]{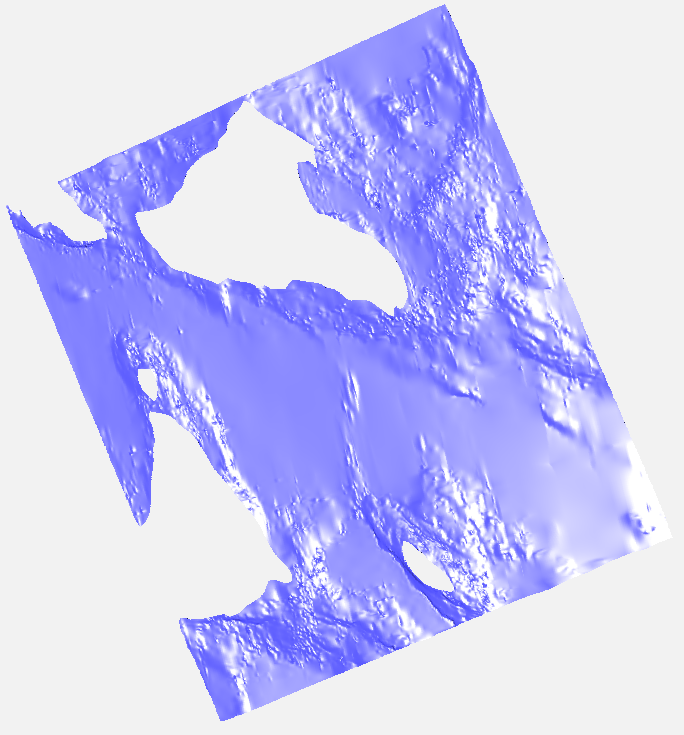} \\
(a) & (b)
\end{tabular}
\caption{\label{fig:sunnmorepts}
  (a) Point set from Sore Sunnm\o re, Norway. (b) Corresponding LR B-spline
  surface. }
\end{figure}
The point cloud, Figure~\ref{fig:sunnmorepts} (a), consists of 11\,150\,110 points obtained from a
881\,100 square meters area at S\o re
Sunnm\o re in Norway. It is a sub sea data set in shallow waters, the depth ranges from
-27.94 meters to -0.55 meters. The point cloud is dense with some holes representing land.
It has some outliers, although less extreme than for the scanline data set. At one position,
two points with the same $xy$-coordinates and differing depth with a distance of 2.38 meters are
given. Thus, the theoretical smallest obtainable maximum distance between the point cloud and the
surface is 1.19 meters. In contrast to the previous
example, the outliers are situated close to other points in the cloud.
The point cloud has a standard deviation of 2.52 meters.

\subsubsection{Result overview}
Despite that it is impossible to reach this resolution, we apply a tolerance of 0.5 meters and run with
a maximum number of iterations of 30 for refinement strategies refining in both parameter
directions and 40 for strategies with alternating parameter directions. The intermediate
stage is where 99.0\% of the points are within the tolerance, meaning that at most 11\,150
points have a larger distance to the surface than 0.5 meters. In
most cases, the iteration runs  the maximum allowed number of steps, however, the improvement in
accuracy in the last steps is small. The differences in the average distances both with regard
to all points and to the out-of-tolerance points are small and the number of new points
being resolved is negligible. A more precise criteria for stopping the iteration should be
applied to avoid unnecessary data size and execution time.

\begin{table}
  \caption{Computation time and number of coefficients. Performance
  ranges and associated strategies for polynomial degrees one to three.}
  \label{tab:Sore_range}
  \begin{tabular}{|l||l|l||l|l||l|l||l|l|}
    \hline
Deg & Min t. & Strategy & Max t. & Strategy & Min cf & Strategy & Max cf & Strategy \\
\hline
1 & 8m47s & eFA & 29m26s & eMc/FB tn & 282\,422 & eFA & 486\,106 & bRB tk \\
\hline
2 & 15m8s & eFA & 1h9m59s & eMcB td & 384\,488 & eFA & 785\,442 & bRB td+k \\ 
\hline
3 & 36m7s & eFA & 3h13m39s & eMlB td & 550\,034 & eFA & 1\,237\,194 & bRB tk \\
\hline
  \end{tabular}
  \end{table}
  Table~\ref{tab:Sore_range} presents the minimum and maximum execution time and
  the lowest and highest number of coefficients in the final surface for polynomial
  degrees one, two and three. The associated refinement strategy is given along 
  with the time and size results. We see that the ranges between best and worst are
  large both with respect to time and size and for all polynomial degrees. The largest
  differences can be found internal to each polynomial degree. 
  Linear polynomials lead in general to lower execution times and smaller surfaces,
  while the cubic polynomials give results at the other end of the scale.

Due to outlier points  neither of the refinement strategies converge with all points within the 
resolution. The minimum number of out-of-tolerance points obtained is 470. We regard all results with
less than 480 points outside the tolerance and a maximum distance of 1.19 meter
as having an accepted accuracy. The least possible maximum distance is reached by all
strategies except bi-quadratic eMcA td (minimum span, alternating parameter directions, threshold
with respect to distance), which finishes with a maximum distance of 1.191 meter.
Several minimum span strategies finish with more than 480 points outside the tolerance. This, in particular,
occurs for strategies with alternating parameter directions, and indicates that the convergence
is slow and the strategy needs more iterations to obtain the expected result.
Minimum span strategies with threshold are also, together with some restricted mesh 
strategies with threshold, the strategies that have the highest computation time and
result in the largest surfaces.

Most strategies continue several iterations after the stage where the number of 
unresolved points and the maximum distance are regarded as satisfactory. The
number of coefficients is increased, in some cases considerably, during these iterations.
The differences are typically largest for the strategies that refine in both
parameter directions and in particular for the restricted mesh strategies. For
instance  the number of coefficients for bRB increase by 122\,418 during the last 13
iterations in the bi-cubic case. 

\subsubsection{Best performances}
\begin{table}
\caption{Sore Sunnm\o re, best results at the intermediate and final stage.}
\label{tab:Sore_summary}
\begin{tabular} {|l?l|l|l|l|l|l|l|}
  \hline
  \multicolumn{8}{|l|}{Intermediate stage} \\
  \hline
  Strategy & Iter & No out & No Coef & Max & Average & Average out & Time\\
\hlinewd{1pt}
\multicolumn{8}{|l|}{Bi-linear} \\
\hline
eMcA tn & 22 & 10\,423 & \cellcolor{Yellow2}201\,216 & 1.679 & 0.063 & 0.594 & 7m32s \\ 
\hline
\cellcolor{Blue2}eFA td & 19 & 11\,093 & 208\,656 & 2.177 & 0.063 & 0.627 & \cellcolor{Yellow2}2m54s \\
\hline
\cellcolor{Blue2}eFA tn & 21 & 9\,967 & 205\,036 & 1.895 & 0.062 & 0.594 & 4m27s\\
\hlinewd{1pt}
\multicolumn{8}{|l|}{Bi-quadratic} \\
\hline
eFA & 19 & 8\,756 & 298\,113 & 2.05 & 0.057 & 0.635 & \cellcolor{Yellow2}5m2s \\
\hline
\cellcolor{Blue2}eFA tn & 20 & 11\,070 & \cellcolor{Yellow2}247\,881 & 1.717 & 0.057 & 0.635 & 5m57s \\
\hlinewd{1pt}
\multicolumn{8}{|l|}{Bi-cubic} \\
\hline
eFB & 10 & 7\,738 & 523\,751 & 1.844 & 0.053 & 0.629 & \cellcolor{Yellow2}11m45s \\
\hline
eMcA tn & 22 & 10\,802 & \cellcolor{Yellow2}357\,114 & 2.673 & 0.045 & 0.609 & 34m51s \\
\hline
\cellcolor{Blue2}eFB tn & 21 & 9\,047 & 375\,301 & 1.799 & 0.055 & 0.6 & 15m22s \\
\hline
\end{tabular}

\vspace*{\floatsep}

\begin{tabular} {|l?l|l|l|l|l|l|l|}
  \hline
  \multicolumn{8}{|l|}{Final stage} \\
  \hline
  Strategy & Iter & No out & No Coef & Max & Average & Time & Tail  \\
\hlinewd{1pt}
\multicolumn{8}{|l|}{Bi-linear} \\
\hline
\cellcolor{Green2}eFA & 40 & 473 & \cellcolor{Yellow2}282\,422 & 1.19 & 0.061 & \cellcolor{Yellow2}8m47s & 21 \\
\hlinewd{1pt}
\multicolumn{8}{|l|}{Bi-quadratic} \\
\hline
\cellcolor{Green2}eFA & 40 & 471 & \cellcolor{Yellow2}384\,488 & 1.19 & 0.056 & \cellcolor{Yellow2}15m8s & 21 \\
\hlinewd{1pt}
\multicolumn{8}{|l|}{Bi-cubic} \\
\hline
\cellcolor{Green2}eFA & 40 & 473 & \cellcolor{Yellow2}550\,034 & 1.19 & 0.053 & \cellcolor{Yellow2}36m7s & 20 \\
 \hline
\end{tabular}
\end{table}
Table~\ref{tab:Sore_summary} presents the best results at the intermediate and final
stage with respect to execution time, surface size and a compound evaluation. In 
contrast to previous test cases do most strategies with the lowest computational time
apply alternating parameter directions, both at the intermediate and final stage. The
least number of coefficients is also obtained with alternating parameter directions.
At the final stage  the full span strategy without threshold delivers the best results
for all degrees and evaluation categories. This result does not change when
observing the situation at the stage where the computation is regarded as converged.
Variations of the full span strategy dominate
also at the intermediate stage. The tail entry at the final stage 
shows that a high number of iterations are executed after the intermediate stage, but the reduction in the number of points
with a larger distance than the tolerance, the reduction
in the maximum distance and the increase in the number of
coefficients are all moderate during these iterations.

\begin{figure}
\centering
\includegraphics[width=\textwidth]{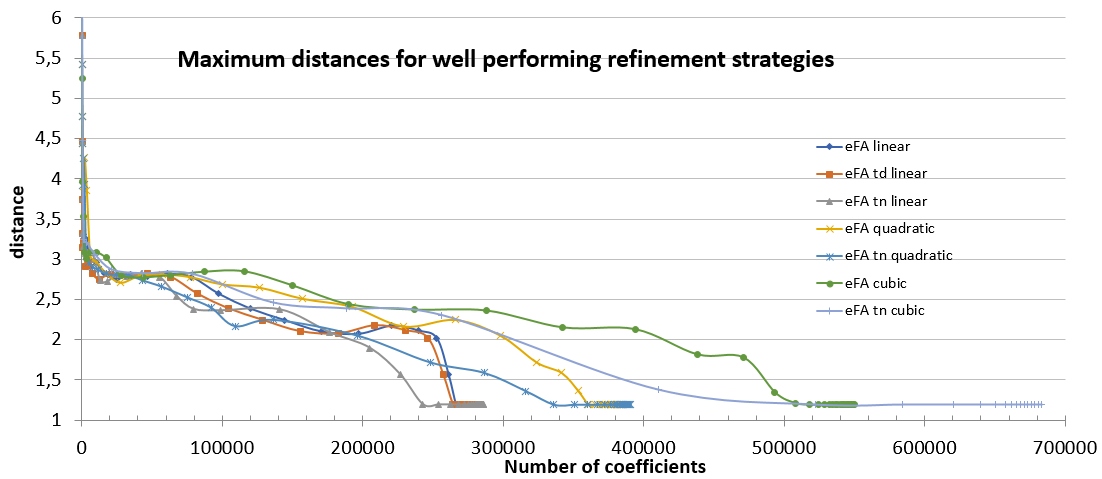}
\caption{
  Maximum distance for selected well performing refinement strategies, all degrees.}
\label{fig:Sore_high_dist}
\end{figure}
%\vspace*{\floatsep}
\begin{figure}
\includegraphics[width=\textwidth]{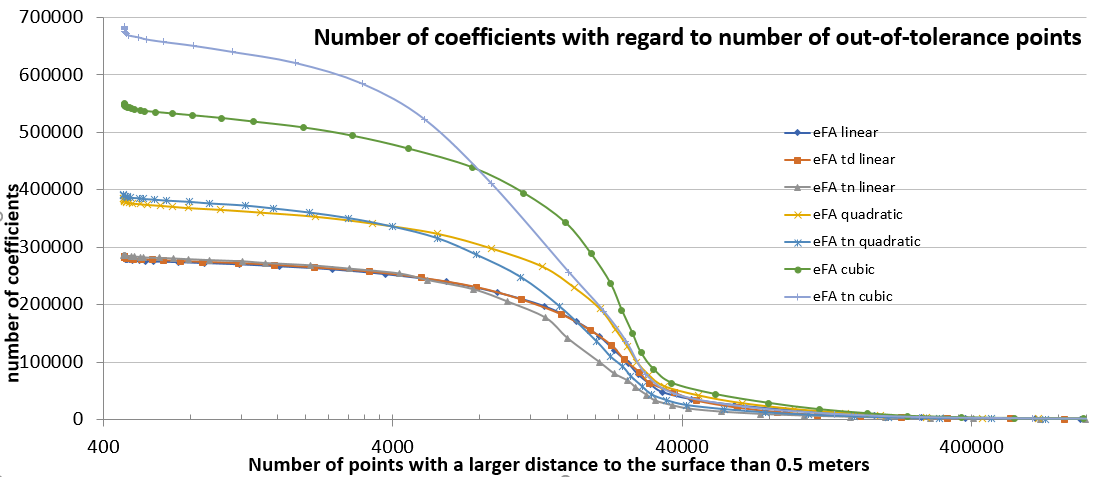}
\caption{Number of coefficients with regard to the number of unresolved points }
\label{fig:Sore_high_out}
\end{figure}
The Figures~\ref{fig:Sore_high_dist} and~\ref{fig:Sore_high_out} give a more
continuous picture of the performance of the strategies selected as best at the
intermediate and final stage. Cubic polynomials have a slower decrease in the 
maximum distance and more coefficients, and this behaviour can be seen 
throughout the approximation process. Cases with linear polynomials are found at the
other end of the specter, but eFA tn (full span, threshold based on number of
out-of-tolerance points), quadratic is compatible with the linear cases
during most of the process.

\subsubsection{Other refinement strategies}
The full span strategies give the best results for this data set as can be seen
from Table~\ref{tab:Sore_summary}.
Restricted mesh refinement strategies behave very similar with and without element 
extension. In both cases  the final number of coefficients and the execution time are
much higher than for corresponding full span strategies. Two restricted mesh strategies without element extension fail to reach the accepted accuracy
in the cubic case. Both apply a 
threshold with respect to distance. All other B-spline based strategies meet this accuracy.  For many cases do minimum span strategies fail
to reach the accepted accuracy. However,  some combinations of minimum span strategies and
threshold have good scores at the intermediate stage. Combinations of a minimum
span strategy and full span strategy turn in several cases  the accuracy result from
not accepted to accepted.  These combinations are only tested with threshold type
tn (based on number of unresolved points). The execution time is, however, considerably higher than
for full span strategies also at the intermediate stage.

\begin{figure}
%\centering
\includegraphics[width=\textwidth]{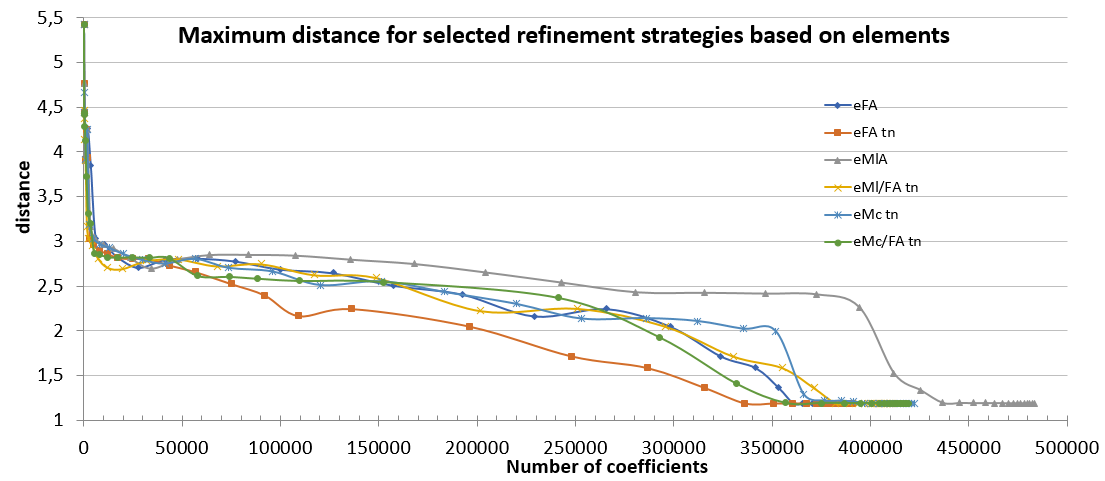}
\caption{Maximum distance for selected alternating refinement strategies.}
\label{fig:Sore_e2_dist}
\end{figure}
%\vspace*{\floatsep}
\begin{figure}
 %   \centering
 \includegraphics[width=\textwidth]{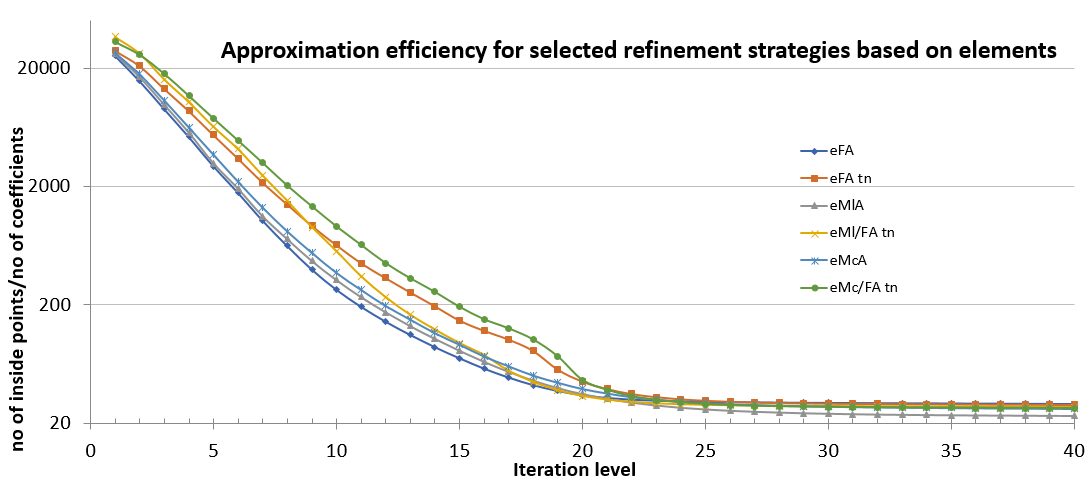}
\caption{Approximation efficiency.}
\label{fig:Sore_e2_eff}
\end{figure}
Figures~\ref{fig:Sore_e2_dist} and~\ref{fig:Sore_e2_eff} compare the performance
of the full span strategy with refinement in alternating parameter directions with
two minimum span strategies. One selects to split the overlapping B-spline with
the largest support (label l) and one combines this criterion with the number of
out-of-tolerance points in the support (label c). Strategies with and without threshold type tn
is shown, and the minimum span strategies are combined with the full span strategy when
the threshold is applied. The maximum distance with respect to the number of
coefficients is shown in Figure~\ref{fig:Sore_e2_dist} and the approximation 
efficiency in Figure\ref{fig:Sore_e2_eff}. Bi-quadratic polynomials are applied.
The maximum distance decreases most rapidly for eFA tn and stays low during the
computation. The approximation efficiency, which is a measure on the number of 
coefficients required for a number of points to be approximated well enough, is also high for eFA tn. The corresponding strategy (full span) without a threshold has a better total score at the
finishing stage, but for most of the computations  the threshold improve the performance.
Also the minimum span strategy eMcFA tn has a higher approximation
efficiency during large parts of the computation. Combining this strategy with
the full span strategy leads to convergence with acceptable time and size figures. 

\subsection{Sea bed, large point cloud}
\subsubsection{Data set}
\begin{figure}
  %\centering
\begin{tabular}{cc}
\includegraphics[width=0.48\textwidth]{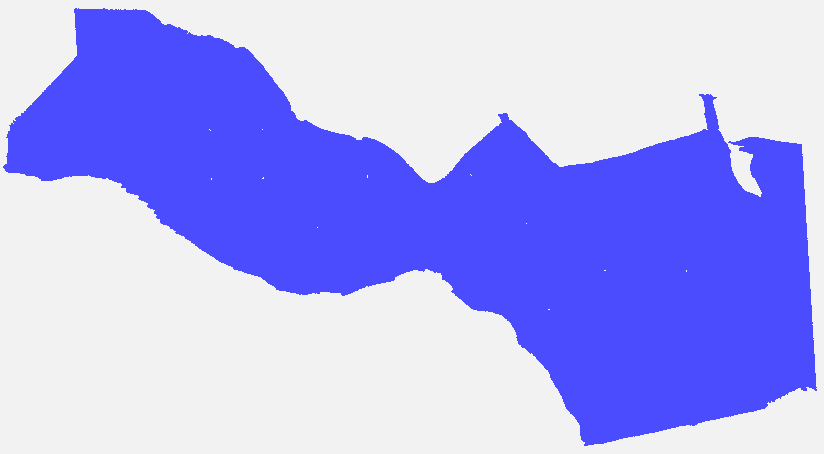}
&\includegraphics[width=0.48\textwidth]{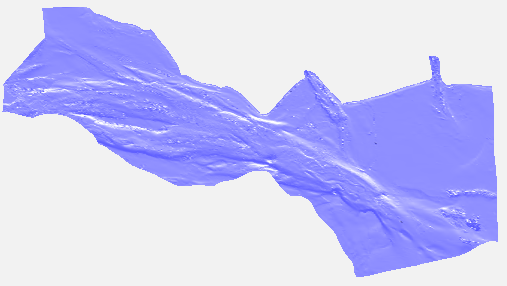} \\
(a) & (b)
\end{tabular}
\caption{\label{fig:largepts}
  (a) Thinned sea bed point cloud. (b) LR B-spline
  surface created from the corresponding full point cloud. }
\end{figure}
The final data set contains 131\,160\,220 well distributed points acquired from the English Channel.
The point cloud is relatively well behaved and represents an area of mostly smooth sea bed,
with some areas of reduced smoothness. Approximation with an LR B-spline surface should, thus,
be well adapted to the properties of the data set. The data set covers about 285 square kilometers of shallow
waters. The elevation ranges from -32.46 meters to 3.26 meters and the standard deviation is 8.04. A thinned version of the point cloud and the surface created with the full span strategy with alternating parameter
directions and threshold
(eFA tn) are shown in Figure~\ref{fig:largepts}.

\subsubsection{Result overview}
Again, we apply a tolerance of 0.5 meters and we let the maximum number of iterations be 40 for refinement
strategies without alternating parameter directions and 60 for alternating directions.
An intermediate stage is defined where {\bf either} 99.9\% of the points are within the
tolerance, which equals 131\,160 points with a larger distance
to the surface than 0.5 meters {\bf or} the maximum distance is less than two meters.

The current point cloud has two complicating factors, the number of data points and the fact that some
points stand out as untypical compared to other points in the neighbourhood. The maximum distance is
for all combinations of refinement strategies and thresholds quickly reduced to about five meters,
then the distance is kept at this level for a number of iterations before it decreases again. Most
strategies converge within the specified number of iterations, a few fail to reach the specified accuracy and some strategies have a very slow
reduction in the number of points outside the tolerance for the last iterations. 

The point cloud 
really distinguishes the strategies both with respect to the final surface size, the execution time and
whether or not the strategy converges. 
\begin{table}
  \caption{Computation time and number of coefficients. Performance
  ranges and associated strategies for polynomial degrees one to three. }
  \label{tab:L_range}
  \begin{tabular}{|l||l|l||l|l||l|l||l|l|}
    \hline
Deg & Min t. & Strategy & Max t. & Strategy & Min cf & Strategy & Max cf & Strategy \\
\hline
\hline
1 & 16m31s & bRB & 1h16m34s & bRA td+k & 235\,033 & bRA & 440\,917 & eMlB td \\
\hline
2 & 30m39s & eFB & 2h31m52s & eMlB td & 291\,753 & bR+eLA td & 729\,432 & eMlB td \\ 
\hline
3 & 56m59s & eFB & 37h44m5s & bRA tk & 374\,490 & eFA & 860\,910 & eMlB td \\
\hline
  \end{tabular}
  \end{table}
  Table~\ref{tab:L_range} shows the ranges in execution time and number of coefficients
  in the final surface for the various polynomial degrees applied. The difference
  between the best and worst result is tremendous especially with regard to time.
In the quadratic case  the  minimum span strategy refining the
largest overlapping B-spline given an element, lead to the surface with the maximum number
of coefficients and the highest execution time. The strategy refines in both parameter direction and
  threshold with respect to the distance between the point cloud and surface is applied (eMlB td). In the linear and
  cubic case  restricted mesh strategies have the poorest
  performance with regard to time.
  The lowest computational time is, for bi-degree two and three, obtained for the full span strategy that
  refines in both parameter directions (eFB).
  \subsubsection{Best performances}
\begin{table}
\caption{Large sea bed point cloud, best results at the intermediate and final stage.}
\label{tab:L_summary}
\begin{tabular} {|l?l|l|l|l|l|l|l|}
  \hline
  \multicolumn{8}{|l|}{Intermediate stage} \\
  \hline
  Strategy & Iter & No out & No Coef & Max & Average & Average out & Time \\ 
\hlinewd{1pt}
\multicolumn{8}{|l|}{Bi-linear} \\
\hline
eFB & 8 & 42\,841 & 223\,855 & 5.412 & 0.075 & 0.678 & \cellcolor{Yellow2}9m45s  \\
\hline
\cellcolor{Blue2}eFB tn & 8 & 88\,791 & 128\,990 & 5.376 & 0.089 & 0.629 & 11m18s  \\
\hline
\cellcolor{Blue2}bRA tk & 16 & 119\,540 & \cellcolor{Yellow2}76\,379 & 5.438 & 0.098 & 0.613 & 17m45s  \\
\hline
\cellcolor{Blue2}bR+eFA tk/n & 16 & 119\,540 & \cellcolor{Yellow2}76\,379 & 5.438 & 0.098 & 0.613 & 17m41s \\
\hlinewd{1pt}
\multicolumn{8}{|l|}{Bi-quadratic} \\
\hline
eFA & 8 & 46\,144 & 274\,288 & 5.392 & 0.066 & 0.690 & \cellcolor{Yellow2}16m8s  \\
\hline
\cellcolor{Blue2}bRB tk & 9 & 130\,965 & \cellcolor{Yellow2}97\,052 & 5.449 & 0.085 & 0.612 & 22m1s \\
\hlinewd{1pt}
\multicolumn{8}{|l|}{Bi-cubic} \\
\hline
eFB & 22 & 71\,175 & 315\,679 & 5.318 & 0.063 & 0.698 & \cellcolor{Yellow2}31m26s \\
\hline
 eMc(+F)A tn & 16 & \cellcolor{Yellow2}126\,520 & 125\,691 & 5.374 & 0.079 & 0.647 & 1h10m34s \\
 \hline
 \cellcolor{Blue2}eFB tn & 8 & 106\,573 & 174\,903 & 5.338 & 0.074 & 0.668 & 34m27s  \\
\hline
\cellcolor{Blue2}eFA tn & 16 & 84\,792 & 145\,020 & 5.482 & 0.074 & 0.663 & 58m53s  \\
\hline
\cellcolor{Blue2}bR+eLB td & 8 & 130\,334 & 145\,641 & 5.488 & 0.074 & 0.655 & 35m11s \\
\hline
\end{tabular}

\vspace*{\floatsep}

\begin{tabular} {|l?l|l|l|l|l|l|l|l|}
  \hline
  \multicolumn{9}{|l|}{Final stage} \\
  \hline
  Strategy & Iter & No out & No Coef & Max & Average & Time & Tail & Oscillations \\
\hlinewd{1pt}
\multicolumn{9}{|l|}{Bi-linear} \\
\hline
bRB & 19 & 0 & 325\,787 & 0.5 & 0.075 & \cellcolor{Yellow2}13m11s & 7 & None \\
  \hline
  \cellcolor{Green2}bRA & 37 & 0 & \cellcolor{Yellow2}235\,033 & 0.5 & 0.082 & 25m11s & 22 & Medium \\
  \hline
  \cellcolor{Green2}bR+eLA & 38 & 0 & 235\,264 & 0.5 & 0.082 & 25m23s & 23 & Medium \\
\hlinewd{1pt}
\multicolumn{9}{|l|}{Bi-quadratic} \\
\hline
eFB & 22 & 0 & 412\,811 & 0.5 & 0.065 & \cellcolor{Yellow2}30m39s & 14 & Low \\
\hline
bR+eLA td & 42 & 0 & \cellcolor{Yellow2}291\,173 & 0.5 & 0.073 & 1h2m7s & 27 & Medium \\
\hline
\cellcolor{Green2}eFA & 36 & 0 & 293\,853 & 0.5 & 0.07 & 44m47s & 21 & Medium \\
 \hlinewd{1pt}
 \multicolumn{9}{|l|}{Bi-cubic} \\
 \hline
 eFB & 21 & 0 & 518\,726 & 0.5 & 0.062 & \cellcolor{Yellow2}56m59s & 13 & Low \\
 \hline
\cellcolor{Green2}eFA & 33 & 0 & \cellcolor{Yellow2}374\,490 & 0.5 & 0.065 & 1h24m20s & 18 & Medium \\
 \hline
\end{tabular}
\end{table}
Table~\ref{tab:L_summary} summarizes the best results at the intermediate and the final stage. 
Different versions of the full span strategy stand out together with some versions of the restricted mesh strategy at
the intermediate stage.
Restricted mesh with threshold tk is chosen as one of the best strategies
at this stage for bi-linear and bi-quadratic polynomials, but these strategies fail to
converge. They are too restrictive in increasing the degrees of freedom in the last part of the
computation. The full span strategy with refinement in
alternating parameter directions (eFA) gets the best score
at the final stage for degree two and three. In the linear case
 the restricted mesh with and without element extension
perform well.

\begin{figure}
\centering
\includegraphics[width=\textwidth]{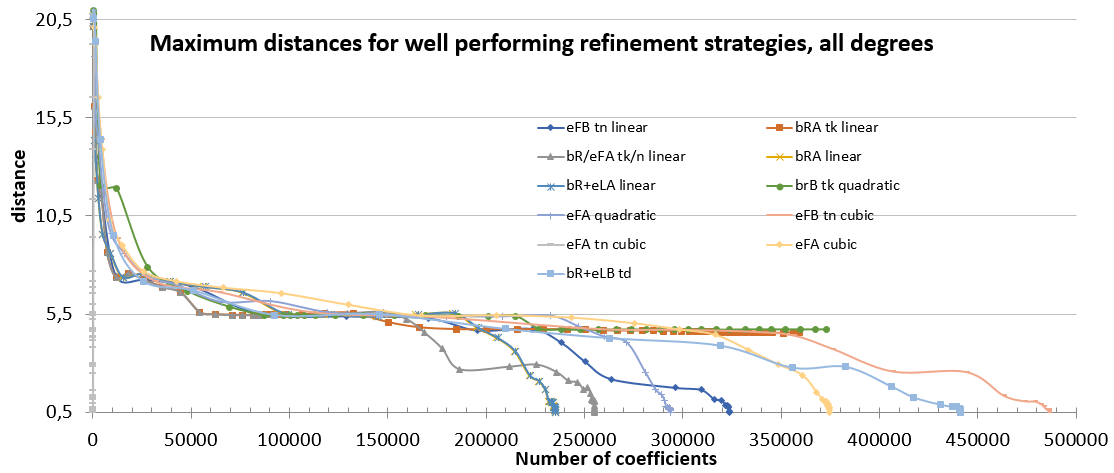}
\caption{Maximum distance for selected well performing refinement strategies, all degrees.}
\label{fig:L_high_dist}
\end{figure}
%\vspace*{\floatsep}
\begin{figure}
\includegraphics[width=\textwidth]{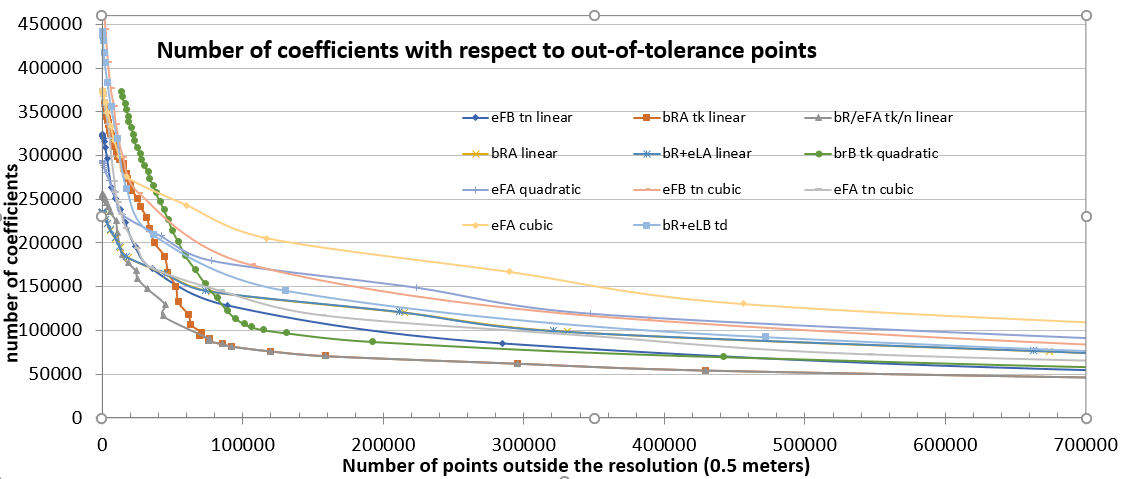}
\caption{Number of coefficients with respect to number of outside points, final stages, all degrees.}
\label{fig:L_high_out}
\end{figure}
Fig.~\ref{fig:L_high_dist} shows the development of the maximum
distance between the point cloud and the surface throughout the computation for the selected methods at the intermediate and
final stage. We see that the pure restricted span strategy with a
threshold (bRA tk) in the bi-linear case and the same strategy with refinement in both parameter directions in the bi-quadratic case reduces the maximum distance
rapidly, but gets stuck at a distance of about 4.5 meters. When  bRA tk linear is combined
with the full span strategy with alternating parameter directions and
threshold (eFA tn),  the computation converges with relatively few coefficients.
The bi-linear surfaces have the least number of coefficients, then come
the bi-quadratic surfaces and finally the bi-cubic ones. The 
difference in the number of coefficients is large even for
selected strategies with the same polynomial degree.

Fig.~\ref{fig:L_high_out} shows the number of coefficients given 
the number of out-of-tolerance points for the same strategies.
Only the last part of the computations are
shown and the development is from right to left. Bi-cubic approximation
leads to more coefficients than the alternative degrees over the
range of out-of-tolerance points. The number of coefficients is
the least in the bi-linear case. The internal ranging of the strategies stay
the same until less than 100\,000 points (less than 0.1\% of the total
number of points) is unresolved. Then there is a shift. This indicates that the strategy to select depends on whether or not the aim is to resolve all
points, something which is not always possible.

\subsubsection{Other refinement strategies}
\begin{table}
\caption{Other selected strategies.}
\label{tab:L_other}
\begin{tabular}{|l?l|l|l|l|l?l|l|l|l|l|}
  \hline
Refine & \multicolumn{5}{l?}{Intermediate stage} & \multicolumn{5}{l|}{Final results} \\
strategy & It & Pt out & No coef  & Max  & Time &  It &
Pt out & No coef  & Max  & Time\\
\hlinewd{1pt}
 \multicolumn{11}{|l|}{bi-linear} \\
 \hline
eFA td & 21 & 107\,275 & 146\,793 & 4.679 &  23m36s & 43 & 0 & 280\,494 & 0.5 & 37m44s \\
\hline
eFA tn & 15 & 123\,366 & 88\,414 & 5.482 & 16m44s & 45 & 0 & 254\,364 & 0.5 & 35m12s \\
\hlinewd{1pt}
 \multicolumn{11}{|l|}{Bi-quadratic} \\
\hline
eFA tn & 15 & 115\,818 & 101\,492 & 5.432 & 30m15s & 46 & 0 & 301\,098 & 0.5 & 1h5m40s \\
\hline
bR+eLA & 15 & 79\,275 & 175\,492 & 5.398 & 27m56s & 41 & 0 & 294\,369 & 0.5 & 50m24s \\
 \hline
bR+eLA td+k & 16 & 83\,659 & 112\,329 & 5.508 & 34m22s & 55 & 0 & 292\,697 & 0.5 & 1h18m9s \\
\hlinewd{1pt}
 \multicolumn{11}{|l|}{bi-cubic} \\
\hline
bRA & 15 & 130\,559 & 193\,111 & 5.438 & 54m41s & 34 & 0 & 387\,211 & 0.5 & 1h32m48s \\
\hline
bR+eLA & 15 & 121\,854 & 197\,582 & 5.438 & 54m7s & 33 & 0 & 386\,891 & 0.5 & 1h28m17s \\
\hline
\end{tabular}
\end{table}
Table~\ref{tab:L_other} presents the figures for some additional strategies at the intermediate and final stages. In the bi-linear case we see that applying a threshold with respect to the number of out-of-tolerances points outperforms the 
distance related threshold when it comes to execution
time and surface size, but the maximum distance at the
intermediate stage is less with threshold of type td.
The full span strategy eFA tn in the bi-quadratic case has less coefficients than eFA, see Table~\ref{tab:L_summary}, at the intermediate stage and more coefficients at the final stage. The execution time is higher when the threshold is applied. The difference in surface size for the element extended restricted mesh strategy with and without a threshold is low at the final stage. 
The increased execution time when the threshold is
applied is more significant. For the bi-cubic polynomials, we see that the effect of adding
element extension to the restricted mesh strategy is small.

Some combinations of strategy and threshold converge for
the minimum span and restricted mesh strategies, but several
fail and the number of
coefficients and the computational time are in general high, in particular
when thresholds are involved. The restricted mesh strategy with threshold tk and td+k struggles most. Switching to the full span strategy helps in some cases, but in general  the swap is done too late. Seen as groups  the full span strategies
and the element extended restricted mesh strategies perform best for
this data set. Thresholds  often reduce the number of coefficients,
but not consistently. The execution time is almost consistently 
increased with threshold.

\section{Summary of results} \label{sec:summary}

This section provides a summary of the results for the individual data sets.
The full span strategy and the structured mesh strategy are very
clearly defined. All other strategies have an element of choice and tuning included in them. Even for the minimum span strategies where
the largest overlapping B-spline for all elements with outside points are subdivided, eMl, there might be more than one candidate B-spline of the same size.

Bi-cubic approximation almost always result in more coefficients
in the resulting surface and higher execution times than corresponding strategies for
lower order polynomials. The pattern regarding 
the relative performance of the strategies for a given data set is quite similar in the
quadratic and the cubic case. The surface size and execution
time is often low in the bi-linear case, but strategies with
some degree of restrictions on the introduction of new knots
are more vulnerable with this degree.

Full span strategies have the most consistent behaviour
for all data sets. The strategy group tends to have
the lowest execution times. The resulting surfaces do
not always have the smallest number of coefficients, but
the results are always among the best.

Minimum span strategies have a high risk of lack of convergence. 
The effect is most pronounced for bi-linear polynomials and reinforced if 
thresholds are applied. Minimum span strategies perform
better in the first part of the computation than for
the final iterations. The combined strategy
for selecting the B-spline to subdivide when an element is
selected perform better than focusing either on the
size of the B-spline support or the number of out-of-tolerance points in the support.

The execution time for the structured mesh strategy is
low, but the surface size is very high for all tested
variations of the strategy, which are omitted for further
testing after two data sets.

The restricted mesh strategy does not have a consistent behaviour for all data sets. It often has a high approximation efficiency in the first and middle part
of the computation, but may need a back-up strategy
based on element to converge completely. Applying a threshold  most often reduces the surface size, but increases the risk of not converging. All variants of the
restricted mesh strategy lead to a high number of
coefficients for the Sore Sunnm\o re data set.

Strategies with alternating parameter directions tend to lead to less coefficients in the final surface. The execution time is most often lowest when refinement is performed in both parameter directions simultaneously,
but the difference is moderate.

Applying a threshold  in general reduces the surface size and increases the computation time, but is not always effective.
A failure to split important B-splines can lead to more coefficients as the
algorithm tries to compensate by performing other and less important splits. The effect of applying a threshold depends on the type. Threshold with respect to distance (td), tends to reduce the maximum distance quite
rapidly, but has little effect on the number of resolved
points. In fact, will strategies with a td threshold often
result in more coefficients than strategies without a threshold. This threshold type  in general have a more
positive effect on the B-spline based strategies than on the element based ones.
A threshold with respect to the number of out-of-tolerance
points (tn and tk) does often lead to good approximation efficiency, that is few coefficients with respect to the
number of resolved points. The effect on the final surface
size varies, the surface can be more lean when applying a threshold, but this is not always the case. It would probably be beneficial if applied thresholds are released
more rapidly than what is done in this study.

Two partly conflicting hypotheses were formulated in 
Section~\ref{sec:refstrat}. Both are confirmed. A reduced
pace in the introduction of new knots results in leaner surfaces,
while a too restrictive approach have the opposite result.

The convergence is typically rapid in the start of the
approximation before it slows down drastically towards the end. The
cost of resolving the last data points is high both with regard
to data size and time consumption.

\section{Conclusion}\label{sec:conclude}
We have performed a study to find out how different strategies for knot selection in adaptive refinement of LR B-spline surfaces, in the context of scattered
data approximation, perform. The tests on the refinement strategies are applied to five different geospatial data sets with a large variety of sizes and
properties. It is not the case that one size fits all. The various strategies
perform differently on different data sets and the performance of each refinement strategy varies throughout the computations. The purpose of approximating a point cloud with an LR B-spline surface is not always to achieve an accurate approximation. 
It can also be to represent the smooth part of the data
points with a smooth surface in order to analyze the
residuals. A very accurate approximation of noisy data is also 
unattractive. The recommendation of the refinement strategy to select depends on the goals for the approximation.

Bi-cubic polynomials should only be used if a surface with
a continuous second derivative is required. A bi-linear surface is a good choice if no smoothness is needed.
We will, in general, recommend the use of bi-quadratic
polynomials.

If an accurate surface is expected should the full span
strategy with alternating parameter directions be applied.
A threshold with respect to the number of points outside the
given tolerance should probably be applied, but there is
room for improvement regarding the actual composition.
Threshold with respect to the maximum distance should not
be applied in this case.

If the aim is to get a smallest possible surface after a
restricted number of iterations or when applying a stop
criterion different from complete convergence, also other
strategies are applicable. The restricted mesh strategy perform well for four out of five test cases. Minimum
span strategies are more seldom selected as the best choice, but are more consistent. In both cases 
a switch to the full span strategy should be applied when the
convergence slows down. The full span strategy is 
applicable also in this case.

Refinement in alternating
parameter directions is recommendable also for early 
finalization of the computation.
Threshold with respect to distance should be used only if
the most significant aim is to reduce the maximum distance
rapidly with less focus on the number of out-of-tolerance points. A threshold with respect to the number of points
is recommended.

The applied study is extensive. More combinations are
possible, but it is questionable whether or not they will provide more insight. The criteria used to define thresholds and restricted strategies have room for improvements.

The study focuses on surface size, execution time and
approximation accuracy. The structure of the polynomial
mesh is not considered although it may have an importance.
B-spline based strategies give a more orderly mesh than
the element based strategies. Furthermore would an analysis on the
size of the accuracy threshold and the best criterion to stop the iteration complement this study in a
good way. Even though a good refinement strategy should not
break down when pushed to the limits does it not necessarily
lead to the best result to insist on an accurate approximation
of all points.

\section*{Acknowledgement}
This work has been supported by the Norwegian research council under grant number 270922.

The data sets are provided by the University of Brest, the UK Hydrographic Office and the Norwegian map authorities.

\end{document}